\documentclass[12pt]{article}
\usepackage{amssymb}
\usepackage{amsfonts}
\usepackage{amsmath}
\usepackage{latexsym}
\usepackage[pdftex]{color}
\usepackage[title]{appendix}
\usepackage{comment}
\usepackage{marginnote} 

\usepackage[a4paper,left=3cm,right=2.5cm,top=3cm,bottom=3cm]{geometry}

\newtheorem{theorem}{Theorem}[section]
\newtheorem{corollary}[theorem]{Corollary}
\newtheorem{lemma}[theorem]{Lemma}
\newtheorem{definition}[theorem]{Definition}
\newtheorem{problem}[theorem]{Problem}
\newtheorem{proposition}[theorem]{Proposition}
\newtheorem{rem}[theorem]{Remark}
\newtheorem{example}[theorem]{Example}

\def\deg{\mathop{\rm deg }\nolimits}
\def\rank{\mathop{\rm rank}\nolimits}
\def\diag{\mathop{\rm diag }\nolimits}

\def\rev{\mathop{\rm rev }\nolimits}

\newcommand{\s}{\ensuremath{\stackrel{s}{\sim}}}

\newcommand{\efe}{\mathbb F}
\newcommand{\FF}{\mathbb F}

\newcommand{\RR}{\mathbb R}

\newcommand{\ba}{{\bf a}}
\newcommand{\bb}{{\bf b}}

\newcommand{\mR}{{\mathcal R}}

\newcommand{\hide}[1]{}

\usepackage{tcolorbox}
\tcbuselibrary{skins}
\tcbuselibrary{breakable}
\tcbset{shield externalize}

\title{Polynomial and rational matrices with \\ the invariant rational
functions and the four \\ sequences of minimal indices prescribed\thanks{Funded by the Agencia Estatal de Investigaci\'on of Spain MCIN/AEI/10.13039/501100011033/ and by ``ERDF A way of making Europe'' of the ``European Union'' through grants PID2021-124827NB-I00, PID2023-147366NB-I00 and RED2022-134176-T,  and by UPV/EHU through grant GIU21/020.}}
\author{Itziar Baraga\~na\thanks{Departamento de Ciencia de la Computaci\'on e I.A., Universidad del Pa\'{\i}s Vasco UPV/EHU, Donostia-San Sebasti\'an, Spain (itziar.baragana@ehu.eus).}, Froil\'{a}n M. Dopico\thanks{Universidad Carlos III de Madrid, ROR: {\tt https://ror.org/03ths8210}, Departamento de Matem\'aticas, Avenida de la Universidad, 30 (edificio Sabatini), 28911 Legan\'es (Madrid), Spain (dopico@math.uc3m.es).}, Silvia Marcaida\thanks{Departamento de Matem\'aticas, Universidad del Pa\'{\i}s Vasco UPV/EHU, Bilbao, Spain (silvia.marcaida@ehu.eus).}, Alicia Roca\thanks{Departamento de Matem\'atica Aplicada, IMM, Universitat Polit\`ecnica de Val\`encia, 46022,Valencia, Spain (aroca@mat.upv.es).}}

\date{\today}

\begin{document}

\maketitle

\begin{abstract}
The complete eigenstructure, or structural data, of a rational matrix $R(s)$ is comprised by its invariant rational functions, both finite and at infinity,  which in turn determine its finite and infinite pole and zero structures, respectively, and by the minimal indices of its left and right null spaces. These quantities arise in many applications and have been thoroughly studied in numerous references. However, other two fundamental subspaces of $R(s)$ in contrast have received much less attention: the column and row spaces, which also have their associated minimal indices. This work solves the problems of finding necessary and sufficient conditions for the existence of rational matrices in two scenarios: (a) when the invariant rational functions and the minimal indices of the column and row spaces are prescribed, and (b) when the complete eigenstructure together with the minimal indices of the column and row spaces are prescribed. The particular, but extremely important, cases of these problems for polynomial matrices are solved first and are the main tool for solving the general problems. The results in this work complete and non-trivially extend  the necessary and sufficient conditions recently obtained for the existence of polynomial and rational matrices when only the complete eigenstructure is prescribed.
\end{abstract}


\medskip

{\bf \small AMS classification codes (2020).} {\small 15A18, 15A21, 15A29, 15A54, 93B10, 93B60}

\section{Introduction}
\textit{ Rational matrices} are matrices whose entries are scalar rational functions in one variable. If all the entries are polynomials, then rational matrices become \textit{ polynomial matrices}. Rational matrices are classical objects in matrix analysis that arise in many applications. For instance, the transfer function matrices of linear multivariable input-output systems are rational matrices \cite{Kail80, Rose70}. Polynomial matrices in their own right also arise in classical applications such as, for example, in the solution of systems of differential equations with constant coefficients \cite{GLR-SIAM-2009}. More recently, rational and polynomial matrices have received a lot of attention in the search for numerical solutions of nonlinear eigenvalue problems (see, for instance, \cite{nlevp-collection, guttel-tisseur-nlep} and the references therein) and in the related context of their linearizations (see, for instance, \cite{local-lin,MMMM-good,MMMM-lin-class} and the references therein).

In addition to their many applications, polynomial and rational matrices have a very rich mathematical theory \cite{GLR-SIAM-2009,Kail80,Rose70} on which there is currently a considerable amount of research activity. A few examples of recent works solving relevant open problems in the area are \cite{amparan2024parametrizing, das-bora, dmytryshyn-geom, root-vectors-vanni-VD}. In this work, we focus on the solution of several inverse problems for polynomial and rational matrices that complete and non-trivially extend  the results in \cite{AmBaMaRo23,AmBaMaRo25,anguasetal2019,DeDoVa15}. These problems are described in the next paragraphs. The definitions of the concepts mentioned in the rest of this introduction can be found in Section \ref{secpreliminaries} and the references therein.

We deal with rational functions in the variable $s$ whose polynomial numerators and denominators have coefficients in a field $\FF$. The set of these rational functions is a field denoted by $\FF (s)$. A rational matrix is \textit{regular} if it is square and its determinant is not identically zero. Otherwise, a rational matrix is said to be \textit{singular}.

The eigenstructure of a regular rational matrix consists of its finite invariant rational functions and its invariant rational functions at infinity. They are revealed by the Smith-McMillan form and the Smith-McMillan form at infinity, respectively, and  fully determine the finite and infinite pole and zero structures of a rational matrix. Polynomial matrices do not have finite poles and their finite zeros are often called eigenvalues in the literature. For polynomial matrices, the invariant rational functions are polynomials and are known as invariant factors, and the structure at infinity is traditionally described by the so-called partial multiplicities of $\infty$. These partial multiplicities of $\infty$ together with the degree give the same information as the invariant rational functions at infinity, but expressed in a different way.

Singular rational matrices have additional structural data related to the fact that they have nontrivial left and/or right null spaces over the field $\FF(s)$. The  rich structures of these spaces are determined by two sequences of nonnegative integers, called right and left minimal indices, respectively, and are important in applications \cite[Section 6.5.4]{Kail80}.
This has motivated to define the \textit{structural data}, sometimes called \textit{complete eigenstructure} or just eigenstructure \cite{vandooren-thesis}, of a rational matrix as the set of  the finite invariant rational functions, the invariant rational functions at infinity, the minimal indices of its right null space and the minimal indices of its left null space. The structural data have been studied in depth in many references, including classical monographs such as \cite{Kail80, Rose70, Vard91}, and, due to their applications, their numerical computation has received considerable attention since the 1970s \cite{vandooren-thesis, vandooren-iee81}.

It must be remarked that the four components of the complete eigenstructure of a rational matrix are not simultaneously revealed by any known canonical form. This makes it difficult to study problems related to it. For instance, the problem of finding a necessary and sufficient condition for the existence of a rational matrix when a complete eigenstructure is prescribed has been recently solved. It was solved in 2015 for polynomial matrices when the field $\FF$ is infinite \cite{DeDoVa15}, in 2019 for general rational matrices with the same restriction over the field \cite{anguasetal2019}, and  in 2024 and 2025 for polynomial and rational matrices over arbitrary fields \cite{AmBaMaRo23, AmBaMaRo25}, and the solutions to these problems are notoriously involved, or rely on results with involved proofs. We must point out, however, that the obtained necessary and sufficient condition is extremely simple and is related to a fundamental result known as the index sum theorem \cite{DeDoMa14, vergheseetal}. See Theorem \ref{theoexistenceDeDoVa152} below and \cite[Theorem 4.1]{anguasetal2019}.

Though the definition of the complete eigenstructure of a rational matrix is well motivated and supported by the relevance in applications of its four components, it can be considered ``unbalanced'' because one of the best known facts in matrix analysis is that any matrix has four fundamental associated subspaces (its left and right null spaces and its column and row spaces \cite{strangfund}). However, the complete eigenstructure only involves the two null spaces. Moreover, the column and row spaces of a rational matrix over the field $\FF (s)$ also have their sequences of minimal indices and, therefore, it is natural to wonder what is the relationship between  the minimal indices of the column space and of the row space with the complete eigenstructure.

Motivated by the discussion above, in this paper we first solve the following problem:
\begin{itemize}
\item[(P1)] find necessary and sufficient conditions for the existence of a rational matrix when the finite invariant rational functions and  the invariant rational functions at infinity are prescribed together with \textit{minimal bases} (and not just minimal indices) of the column and row spaces.
\end{itemize}
Then, on the basis of the solution of (P1) and of other results available in the literature, we solve the following problems:
\begin{itemize}
    \item[(P2)] find necessary and sufficient conditions for the existence of a rational matrix when the finite invariant rational functions, the invariant rational functions at infinity, the minimal indices of the column space and the minimal indices of the row space are prescribed.
    \item[(P3)] find necessary and sufficient conditions for the existence of a rational matrix when a complete eigenstructure and the minimal indices of the column space and the minimal indices of the row space are prescribed.
\end{itemize}
Problem (P2) can be seen as the counterpart of the problem solved in \cite{AmBaMaRo23, anguasetal2019, DeDoVa15} mentioned above when the minimal indices of the left and right null spaces are replaced by the minimal indices of the column and row spaces in the prescribed data. Thus, both problems involve four lists of prescribed data, whilst Problem (P3) involves the minimal indices of the column and row spaces  in addition to the complete eigenstructure, and thus a total of six lists of prescribed data, which is a considerable amount of prescribed information.
 However, it should be noted that rational matrices have a very rich structure and that these six lists are not the only data of interest associated with them. For instance, other indices of interest are the right and left Wiener--Hopf factorization indices at infinity \cite{FuWi79}, which here are not going to be taken into account.

The solutions to the problems posed above are achieved in two steps. First, the problems are solved for polynomial matrices, paying particular attention to the prescription of the degree, and, then, are solved for general rational matrices using the solutions of the polynomial cases. We remark that the obtained necessary and sufficient conditions are rather simple, though not as simple as the constraint related to the index sum theorem (which, as we will see, still plays a role), that a majorization relation arises among the conditions, and that the sufficiency of such conditions only holds when the field $\FF$ is algebraically closed. An example will illustrate that such restriction over the field is unavoidable.

The column and row spaces, and their minimal indices, of rational matrices have received much less attention in the literature than the left and right null spaces and their associated minimal indices. Despite this, they have appeared in some interesting problems. For example, given a rational matrix $G(s)$, the minimal indices of the column space of the compound matrix $\left[ \begin{smallmatrix} I \\ G(s) \end{smallmatrix}\right]$ are important in the study of controller canonical realizations  of transfer function matrices \cite[Section 7]{Fo75}. In a different context, the minimal bases of the column and row spaces of a polynomial matrix play a fundamental role in the so-called minimal rank factorizations recently introduced in \cite{DmDoVa23}.

The rest of the paper is organized as follows. In Section \ref{secpreliminaries} we summarize previously known concepts and results needed throughout the paper and introduce the corresponding notations. In Section \ref{problems} complete and rigorous formulations of the problems (P1), (P2) and (P3) are introduced, considering for brevity only the case of polynomial matrices. In Section \ref{sec.factorization}, some results on factorizations of polynomial matrices are developed and are used to reformulate problem (P1) in a more appropriate form to achieve a solution. Solutions to the problems (P1), (P2) and (P3) for polynomial matrices are given in Theorems \ref{theomainlk2}, \ref{cormainlk} and \ref{cormainlkcv} in Section \ref{main}, respectively. In Section \ref{sec.rational} solutions to the three problems for general rational matrices are presented in Theorems \ref{theomainlk2_rat}, \ref{cormainlk_rat} and \ref{cormainlkcv_rat}. The conclusions and some open problems are discussed in Section \ref{sec.conclusions}. Finally, the long proof of a technical determinantal lemma needed in the proof of one of the main results is presented in Appendix \ref{secappend}.

\section{Preliminaries}\label{secpreliminaries}
Throughout this work $\FF$ denotes a field and $\overline{\FF}$ its algebraic closure. Unless otherwise indicated, the field $\FF$ is arbitrary. $\FF[s]$ stands for the ring of univariate polynomials in the variable $s$ with coefficients in $\FF$ and  $\FF(s)$ for the field of fractions of $\FF[s]$. For two scalar polynomials $\alpha_1 (s), \alpha_2 (s) \in \FF [s]$,  $\alpha_1 (s) \mid \alpha_2 (s)$ denotes that $\alpha_1 (s)$ divides $\alpha_2 (s)$. The elements of $\FF(s)$ with the degree of the numerator at most the degree of the denominator are called \textit{proper rational functions}, and $\FF_{pr}(s)$ denotes the ring of proper rational functions over $\efe$; if the degree of the numerator is strictly less than that of the denominator they are called \textit{strictly proper rational functions}. Fractions in $\FF_{pr}(s)$ with numerator and  denominator of equal degree are called \textit{biproper rational functions} and are the units of $\FF_{pr}(s)$.

We denote by  $\efe(s)^p$ the set of column vectors with $p$ components in $\efe(s)$, and by $\efe^{p\times m}$ ($\efe[s]^{p\times m}$, $\efe(s)^{p\times m}$, $\efe_{pr}(s)^{p\times m}$) the set of $p\times m$ matrices with entries in $\efe$ (respectively $\efe[s]$, $\efe(s)$, $\efe_{pr}(s)$). Matrices in $\efe[s]^{p\times m}$ are called \textit{polynomial matrices} or \textit{matrix polynomials} indistinctly. The \textit{degree} of a polynomial matrix, denoted $\deg(\cdot)$, is the highest degree of its entries.  We define the degree of the zero polynomial to be $-\infty$ and take $s^{-\infty} \equiv 0$. A non singular   polynomial matrix is \textit{unimodular}  if its inverse is also polynomial.
Equivalently, unimodular matrices are square polynomial matrices whose determinants are non zero elements of $\FF$. Matrices  in $\efe(s)^{p\times m}$ are known as \textit{rational matrices}, and
matrices  with entries in $\efe_{pr}(s)$ are termed as \textit{proper rational matrices}.  Invertible matrices in $\efe_{pr}(s)^{p\times p}$, that is, non singular proper rational matrices whose inverses are also proper, are called \textit{biproper}.  Equivalently, biproper matrices are square proper rational matrices whose determinants are non zero biproper rational functions. The matrices in $\efe_{pr}(s)^{p\times m}$ whose entries are strictly proper rational
functions are called \textit{strictly proper rational matrices}.

The \textit{rank} of a rational matrix $R(s)$ over the field $\FF (s)$, denoted $\rank(R(s))$, is equal to the largest order of the non identically zero minors of $R(s)$. It is also  known in the literature as normal rank. A rational matrix $R(s)$ can be uniquely decomposed as a sum $R(s) = P(s) + R_{sp} (s)$, where $P(s)$ is a polynomial matrix and $R_{sp} (s)$ is a strictly proper rational matrix. $P(s)$ and $R_{sp} (s)$ are called the polynomial and the strictly proper parts of $R(s)$, respectively.

Since polynomial matrices are also rational matrices, concepts defined for rational matrices can be applied to polynomial matrices. However, the literature on polynomial matrices often uses a somewhat different nomenclature and slightly different definitions for some notions, for instance the infinite structure. Therefore,  in this section we will pay particular attention to polynomial matrices.
This is also connected to the fact that the results of this paper are obtained first for polynomial matrices and then  extended to rational matrices.

Along the paper we use the notion of column proper or column reduced polynomial matrix (\cite[Section 6.3.2] {Kail80}; see also  \cite[Definition 1.10]{Vard91}, \cite[Definition 2.5.6]{Wo74}), which is recalled in this paragraph. Let $P(s)\in \FF[s]^{m \times n}$ be a polynomial matrix. For $1\leq i \leq n$, we denote by $d_i$ the degree of the $i$-th column of $P(s)$, and $d_1, \ldots , d_n$ are called the column degrees of $P(s)$. The matrix $P(s)$ can  be written as
\begin{equation*}\label{reduced}
P(s)=P_h \diag(s^{d_1}, \ldots, s^{d_n})+ L(s),
\end{equation*}
where  $P_h \in \FF^{m\times n}$ is the  highest column degree coefficient matrix of $P(s)$, i.e., a matrix whose $i$-th column comprises the coefficients of $s^{d_i}$ in the $i$-th column of $P(s)$, and $L(s) \in \FF[s]^{m\times n}$ is a polynomial matrix collecting the remaining terms. For the non zero columns of $P(s)$, $L(s)$ has lower column degrees than the corresponding ones in $P(s)$. The  matrix $P(s)$ is called  \textit{column proper} or \textit{column reduced} if
 $P_h$ has full rank, i.e., $\rank (P_h) =\min\{m, n\}$.
  In the case that $m=n$, note that if $P(s)$ is column proper, then $\deg(\det(P(s)))=\sum_{i=1}^n d_i$. A polynomial matrix is called  \textit{row proper} or \textit{row reduced} if its transpose is column reduced.

In this work we frequently deal with sequences of integers. Let $a_1, \dots,  a_m$ be a  sequence of integers. Whenever we write  $\ba=(a_1, \dots,  a_m)$, we  understand that   $a_1\geq \dots \geq a_m$. If $a_m\geq 0$, then
$\ba=(a_1, \dots,  a_m)$ is  called a \textit{ partition}.

Let   $\ba= (a_1,  \ldots, a_m)$ and $\bb= (b_1, \ldots, b_m)$ be  two sequences of integers.  Following \cite{MOA11}, it is said  that   $\ba$ is \textit{majorized} by $\bb$ (denoted by $\ba \prec \bb$) if $\sum_{i=1}^k a_i \leq \sum_{i=1}^k b_i $ for $1 \leq k \leq m-1$ and $\sum_{i=1}^m a_i =\sum_{i=1}^m b_i$.

Given a sequence of integers $a_1, \dots,  a_m$, we agree that $\sum_{i=p}^q a_i=0$ when $p>q$.

\subsection{Equivalence relations on rational matrices}

Two rational matrices $R_1(s), R_2(s)\in \FF(s)^{m\times n}$ are
 \textit{unimodularly equivalent} if there exist unimodular matrices  $U(s)\in \FF[s]^{m\times m}$ and $V(s)\in \FF[s]^{n\times n}$ such that
$R_2(s)=U(s)R_1(s)V(s)$.
A rational matrix  $R(s)\in \FF(s)^{m\times n}$ of $\rank(R(s))=r$   is unimodulary equivalent to its \textit{Smith--McMillan form} \cite[p. 109]{Rose70}
\begin{equation} \label{eq.smithmcmillan}
M(s)=
\begin{bmatrix}\diag\left(\frac{\epsilon_1(s)}{\psi_1(s)},\dots,\frac{\epsilon_r(s)}{\psi_r(s)}\right)&0\\0&0\end{bmatrix},
\end{equation}
where $\epsilon_1(s)\mid  \dots \mid \epsilon_r(s)$ and $\psi_r(s)\mid\dots\mid\psi_1(s)$ are monic polynomials, and $\frac{\epsilon_1(s)}{\psi_1(s)},\dots,\frac{\epsilon_r(s)}{\psi_r(s)}$ are irreducible rational functions called the \textit{(finite) invariant rational functions} of $R(s)$. The Smith-McMillan form is a canonical form for the unimodular equivalence.

Two rational matrices $R_1(s), R_2(s)\in \FF(s)^{m\times n}$ are
 \textit{equivalent at  infinity} if there exist biproper  matrices $B_1(s)\in\FF_{pr}(s)^{m\times m}$ and $B_2(s)\in\FF_{pr}(s)^{n\times n}$
such that
$R_2(s)=B_1(s)R_1(s)B_2(s)$.
 A canonical form for the equivalence at infinity of rational matrices is the  \textit{Smith-McMillan  form at infinity} (see \cite[Section 5]{AmMaZa15} or \cite[Theorem 3.13]{Vard91}). That is, a rational matrix
$R(s)\in\FF(s)^{m\times n}$ of $\rank (R(s))=r$  is equivalent at infinity to
\begin{equation} \label{eq.infsmithmcmillan}
\begin{bmatrix}\diag\left(s^{-q_1},\ldots,s^{-q_r}\right)&0\\0&0\end{bmatrix},
\end{equation}
where  $q_1\leq\cdots\leq q_r$ are integers called the \textit{invariant orders at infinity} of $R(s)$, and $s^{-q_1},\ldots,s^{-q_r}$ are the \textit{invariant rational functions at infinity} of $R(s)$. Note that the smallest invariant order at infinity $q_1$ is minus the degree of the polynomial part of $R(s)$ if this polynomial part is non zero, and $q_1$ is positive otherwise \cite[p. 102]{Vard91}.

If $P(s)\in \FF[s]^{m\times n}$ is a polynomial matrix of $\rank (P(s))=r$, the denominators in \eqref{eq.smithmcmillan} are $\psi_1(s)=\dots=\psi_r(s)=1$, that is, its invariant rational functions are just the numerators $\epsilon_1(s)\mid  \dots \mid \epsilon_r(s)$, called the \textit{invariant factors} of $P(s)$, and $M(s)$ reduces to its \textit{Smith form} \cite[Section 7.5]{LaTi85}.

Let $P(s)=P_ds^d+\dots+P_1s+P_0 \in \FF[s]^{m\times n}$ be a polynomial matrix, where $P_i \in \FF^{m\times n}$ for $i= 0, \ldots , d$ and $P_d \ne 0$. It is said that $\lambda \in \overline{\FF}\cup \{\infty\}$ is an \textit{eigenvalue} of $P(s)$ if $\rank(P(\lambda))< \rank(P(s))$, where we understand that $P(\infty)=P_d$. The set of eigenvalues of $P(s)$ is denoted by $\Lambda(P(s))$. We define the \textit{reversal polynomial matrix} of $P(s)$ as
$\rev P(t)=t^d P(1/t)$. It turns out that $\infty$ is an eigenvalue of $P(s)$
if and only if $0$ is an eigenvalue of $\rev P(t)$, because $\rev P (0) = P_d$.

Let $\alpha_1(s) \mid \cdots \mid \alpha_r(s)$ be the invariant factors of the polynomial matrix $P(s) \in \FF[s]^{m\times n}$. Recall that (see, for instance, \cite[p. 261]{LaTi85})
\begin{equation}\label{minork}
\alpha_1(s)\cdots\alpha_k(s)=\gcd\{m_k(s):  \mbox{$m_k(s) = $ minor of order  $k$ of $P(s)$}\},\ 1\leq k\leq r,
\end{equation}
where $\gcd\{\cdot\}$ stands for monic greatest common divisor.
Factorizing the invariant factors of $P(s)$ as products of irreducible polynomials over $\overline{\FF}$,  we can write
$$
\alpha_{i}(s)=\prod_{\lambda\in \Lambda(P(s))\setminus\{\infty\}}(s-\lambda)^{n_i(\lambda, P(s))}, \quad 1\leq i \leq r.
$$

The factors  $(s-\lambda)^{n_i(\lambda, P(s))}$ with $n_i(\lambda, P(s))>0$ are the \textit{elementary divisors}  of $P(s)$ over $\overline{\FF}$ corresponding to $\lambda$, and the non negative integers
$
n_1(\lambda, P(s))\leq \dots \leq n_r(\lambda, P(s))
$
are the \textit{partial multiplicities} of $\lambda$ in $P(s)$.

The  \textit{infinite elementary divisors}  of $P(s)$ and the \textit{partial multiplicities of $\infty$} in $P(s)$ are the elementary divisors  of $\rev P (t)$ corresponding to $0$ and the partial multiplicities of $0$ in $\rev P(t)$, respectively \cite[pp. 184-185]{GLR-SIAM-2009}. For simplicity, we will denote $f_i=n_i(0, \rev P(t))$ for $1\leq i \leq r$.

Whenever $\lambda\in\overline{\FF}\cup\{\infty\}$ is not an eigenvalue of $P(s)$ we define the partial multiplicities of $\lambda$ as $n_1(\lambda, P(s)) = \cdots = n_r(\lambda, P(s)) = 0$.

For $\lambda\in\overline{\FF}$, the polynomials $(s-\lambda)^{n_i(\lambda, P(s))}$ with $n_i(\lambda, P(s))\geq 0$  are  called  the \textit{local invariant rational functions at $s-\lambda$ of the polynomial matrix $P(s)$}  \cite[Section 4]{AmMaZa15}. Observe that $P(s)$ has $r = \rank (P(s))$ local invariant rational functions for each $\lambda\in\overline{\FF}$, with $(s-\lambda)^{n_i(\lambda, P(s))} = 1$ whenever $n_i(\lambda, P(s)) =0$. However, the elementary divisors corresponding to an eigenvalue $\lambda \in \overline{\FF}$ are always different from $1$ and, so, there may be less than $r$.

The above definition can be extended to any rational matrix $R(s) \in \FF (s)^{m\times n}$ as follows. Let $\Omega (R(s))$ be the set of roots in $\overline{\FF}$ of the numerators and the denominators of the invariant rational functions of $R(s)$ in \eqref{eq.smithmcmillan}. Then, we can write
$$
\frac{\epsilon_i (s)}{\psi_i (s)} = \prod_{\lambda\in \Omega(R(s))} (s-\lambda)^{n_i(\lambda, R(s))}, \quad 1\leq i \leq r,
$$
where the integers
$
n_1(\lambda, R(s))\leq \dots \leq n_r(\lambda, R(s))
$
can be negative, zero or positive. Moreover, we define $n_1(\lambda, R(s)) = \cdots = n_r(\lambda, R(s)) = 0$ for $\lambda \notin \Omega (R(s))$.  Thus, for any $\lambda\in\overline{\FF}$, the rational functions $(s-\lambda)^{n_i(\lambda, R(s))}$ are defined and are called the  \textit{local invariant rational functions at $s-\lambda$ of $R(s)$}. See \cite[Section 4]{AmMaZa15} for more details.

It is known \cite[Proposition 6.14]{AmMaZa15} that the partial multiplicities of $\infty$, $f_i$, the invariant orders at infinity, $q_i$, and the degree $d$ of a polynomial matrix $P(s)$ with $r = \rank (P(s))$ are related as follows
\begin{equation}\label{eqqf}
  f_i=q_i+d,\quad 1\leq i \leq r.
  \end{equation}
Observe that the smallest partial multiplicity of $\infty$ in a polynomial matrix is always $f_1 = 0$ because $q_1=-d$. The equalities in \eqref{eqqf} imply that two polynomial matrices are equivalent at infinity if and only if they have the same degree and the same partial multiplicities of $\infty$, since these quantities determine the invariant orders at infinity and vice versa. We emphasize that the condition $f_1 = 0$ must appear in all the results of this work that involve the partial multiplicities of $\infty$ of a polynomial matrix.

For $\FF=\RR$, it was proved in \cite[p. 102]{Vard91} that the invariant orders at infinity $q_1 \leq \cdots \leq q_r$ of a polynomial matrix $P(s)$ satisfy
\begin{equation}\label{eqminors}
\sum_{i=1}^kq_i=-M_k(P(s)), \quad 1\leq k\leq r,
\end{equation}
where
\begin{equation} \label{def.MkP}
M_k(P(s))=
\max\{\deg(m_k(s)): \mbox{$m_k(s)= $ minor of order  $k$ of $P(s)$}\}.
\end{equation}
This result can be trivially extended to any field.
From (\ref{eqqf})  and (\ref{eqminors}) we obtain
\begin{equation}\label{eqmaxminP2}
M_k(P(s))=kd-\sum_{i=1}^kf_i, \quad 1 \leq k \leq r.
 \end{equation}

\subsection{Rational subspaces and minimal indices associated to a rational matrix}

The definitions of the four sequences of minimal indices of a rational matrix are based on the notions of minimal bases and indices of a rational subspace introduced in \cite{Fo75}. We recall first these notions.  A subspace of  $\FF(s)^p$ over the field of rational functions $\FF (s)$ is called a rational subspace. For any rational subspace $\mathcal{V}$ of $\FF(s)^p$, it is possible to find a \textit{polynomial basis}, i.e., a basis consisting of polynomial vectors; it is enough to take an arbitrary basis and multiply each vector
by the  least common multiple of the denominators of its  entries.
The \textit{order of a polynomial basis} of $\mathcal{V}$ is the sum of
the degrees of its vectors. A \textit{minimal basis} of  $\mathcal{V}$
is  a polynomial basis with least order among the
polynomial bases of $\mathcal{V}$. The decreasingly ordered list of degrees of the polynomial vectors of any minimal basis  of $\mathcal{V}$ is always the same \cite[Main Theorem p. 495]{Fo75}. These degrees are  called the \textit{minimal indices} of $\mathcal{V}$.

We present now the null-space or singular structure of a rational matrix \cite[Section 6.5.4]{Kail80}. Let $R(s)\in\FF(s)^{m\times n}$ be a rational matrix. Denote by $\mathcal{N}_r (R(s))$ and $\mathcal{N}_\ell (R(s))$ the rational \textit{right and left null-spaces} of $R(s)$, respectively, i.e.,
\begin{align*}
\mathcal{N}_r (R(s)) & =\{x(s)\in\FF(s)^{n}: R(s)x(s)=0\}, \\
\mathcal{N}_\ell (R(s))& =\{x(s)\in\FF(s)^{m}: x(s)^TR(s)=0\}.
\end{align*}
A \textit{right (left) minimal basis} of $R(s)$ is a minimal basis of $\mathcal{N}_r (R(s))$ ($\mathcal{N}_\ell (R(s))$).
The \textit{right  (left) minimal indices} of $R(s)$ are the minimal indices of
$\mathcal{N}_r (R(s))$ ($\mathcal{N}_\ell (R(s))$). We will work with the right (left) minimal indices ordered decreasingly. Notice that a rational matrix $R(s)\in \FF(s)^{m\times n}$ of $\rank(R(s))=r$ has $n-r$ right minimal indices and $m-r$ left minimal indices, by the rank-nullity theorem.

The sequences of right and left minimal indices of a rational matrix $R(s)$ form the \textit{null-space or singular structure} of $R(s)$. This structure has applications in different problems arising in linear systems theory \cite[Section 6.5.4]{Kail80}.

The right (left) minimal indices of $R(s)$ are also called in the literature column (row) minimal indices of $R(s)$ (see, for instance, \cite{AmBaMaRo23}). We do not use this nomenclature in this paper to avoid possible confusions with the minimal indices defined below.

Given a rational matrix
$R(s)\in\FF(s)^{m\times n}$ of rank  $r$
we  denote by  $\mR(R(s))$ ($\mR(R(s)^T)$)
\textit{the rational subspace} of $\FF(s)^{m}$ ($\FF(s)^{n}$)
\textit{spanned by the columns (rows)} of $R(s)$, which is called the \textit{column (row) space} of $R(s)$. That is,
\begin{align*}
\mR(R(s)) & = \{ R(s) x(s) : x(s) \in \FF(s)^n\}, \\
\mR(R(s)^T) & = \{ R(s)^T x(s) : x(s) \in \FF(s)^m\}.
\end{align*}
The \textit{col-span (row-span) minimal indices of $R(s)$} are the minimal indices of $\mR(R(s))$ ($\mR(R(s)^T)$). Note that $R(s)$ has $r$ col-span minimal indices and $r$ row-span minimal indices.  We will work with the col-span (row-span) minimal indices ordered decreasingly.

In contrast with the minimal indices and bases of $\mathcal{N}_r (R(s))$ and $\mathcal{N}_\ell (R(s))$, the minimal indices and bases of the column and row spaces of rational matrices have not received much attention in the literature, as far as we know. An exception is the recent paper \cite{DmDoVa23}, where they have played a fundamental role
for constructing rank revealing factorizations of polynomial matrices with
good properties.

\subsection{Polynomial and rational matrices with prescribed eigenstructure}
\label{sec.previousprescribed}

The \textit{eigenstructure} of a polynomial matrix is formed by the invariant factors, the partial multiplicities of $\infty$, and the right and left minimal indices. The \textit{eigenstructure} of a rational matrix is formed by the finite invariant rational functions, the invariant rational functions at infinity, and the right and left minimal indices.  This name was introduced in \cite{vandooren-thesis} and has been used in different papers as \cite{DeDoVa15,vd-dewilde83}. Knowing the eigenstructure of a polynomial matrix is equivalent to knowing its eigenstructure as a rational matrix (see the comments after (\ref{eqqf})). In \cite{AmBaMaRo23},  a characterization of the existence of a polynomial matrix with prescribed eigenstructure is provided over arbitrary fields, generalizing the result proved over infinite fields in  \cite{DeDoVa15}. This result is presented in the next theorem.

\begin{theorem}[{\rm \cite[Theorem 3.1]{AmBaMaRo23}, \cite[Theorem 3.3]{DeDoVa15} for infinite fields}]\label{theoexistenceDeDoVa152}
Let $m$, $n$, $r\leq \min\{m,n \}$ be positive integers and $d$ a non negative integer.
Let
$\alpha_1(s)\mid \dots \mid \alpha_r(s)$ be monic polynomials in $\FF[s]$.
Let $(f_r, \ldots, f_{1})$, $(d_1, \ldots, d_{n-r})$, $(v_1, \ldots, v_{m-r})$
be partitions of non negative integers.
Then, there exists $P(s)\in \FF[s]^{m\times n}$, $\rank(P(s))= r$,
$\deg(P(s))= d$,  with invariant factors $\alpha_1(s), \dots, \alpha_r(s)$, $f_1, \dots, f_r$ as partial multiplicities of $\infty$, and right and left minimal indices
$d_1, \dots, d_{n-r}$ and $v_1, \dots, v_{m-r}$, respectively, if and only if $f_1=0$ and
\begin{equation}\label{eqIST}
\sum_{i=1}^{n-r}d_i+\sum_{i=1}^{m-r}v_i+\sum_{i=1}^{r}f_i+\sum_{i=1}^{r}\deg(\alpha_i)=rd.
\end{equation}
\end{theorem}

An analogous result for the existence of a rational matrix with prescribed  eigenstructure was proved in \cite[Theorem 4.1]{anguasetal2019} over infinite fields, and in \cite[Theorem 2.4]{AmBaMaRo25} for arbitrary fields.

Note that the necessity of condition \eqref{eqIST} in Theorem \ref{theoexistenceDeDoVa152} is the well-known \textit{ index sum theorem} (see \cite[Theorem 6.5]{DeDoMa14} or \cite[Theorem 3]{vergheseetal}). This motivates that condition \eqref{eqIST}  is also known as the index sum theorem constraint \cite{anguasetal2019}.

This paper has two principal aims, both related to Theorem \ref{theoexistenceDeDoVa152} and \cite[Theorem 4.1]{anguasetal2019}. First, to obtain characterizations of the existence of polynomial and rational matrices as in Theorem \ref{theoexistenceDeDoVa152},  when the col-span and row-span minimal indices are prescribed instead of the right and left minimal indices. Second, to obtain characterizations of the existence of polynomial and rational matrices when the col-span and row-span minimal indices are prescribed in addition to the right and left minimal indices.

\subsection{Properties of minimal bases and indices} In this paper, minimal bases of rational subspaces are often arranged as columns of matrices. Therefore, for brevity, we abuse of the terminology and say that a polynomial matrix $N(s)\in\FF[s]^{m\times r}$  is a \textit{minimal basis}  if its columns form a minimal basis of $\mR(N(s))$. A characterization of a minimal basis particularly useful in this work is recalled in  the next theorem.

\begin{theorem}[{\rm \cite[Main Theorem-2. p. 495]{Fo75}}]\label{mincharact}
A polynomial matrix $N(s)\in\FF[s]^{m\times r}$ is
a minimal basis  if and only if
$N(s_0)$ has full column rank for all $s_0\in\overline{\FF}$ and $N(s)$ is column reduced.
\end{theorem}

\begin{rem}\label{remmincharact}
{\rm Observe that $N(s)\in\FF[s]^{m\times r}$ satisfies that
  $N(s_0)$ has full column rank for all $s_0\in\overline{\FF}$ if and only if
   $N(s)$ has $r$ invariant factors equal to $1$ or, equivalently, if and only if the Smith form of $N(s)$ is $\begin{bmatrix}I_r\\0\end{bmatrix}$, and this condition implies that $N(s)\in\FF[s]^{m\times r}$ has full column rank. In the square case $m=r$, $N(s_0)$ is non singular for all $s_0\in\overline{\FF}$ if and only if $N(s)$ is unimodular. \hfill $\Box$
   }
\end{rem}

The next result relates two minimal bases of the same rational subspace.
\begin{proposition}\label{eqbasis}
  Let $N_1(s)$, $N_2(s)\in\FF[s]^{m\times r}$ be minimal bases. Then \allowbreak
  $\mR(N_1(s))  = \mR(N_2(s))$
  if and only if there exists a unimodular matrix $U(s)\in\FF[s]^{r\times r}$ such that
  $N_1(s)U(s)=N_2(s)$.
\end{proposition}

{\bf Proof.} If $\mR(N_1(s))= \mR(N_2(s))$, there exists $G(s)\in\FF(s)^{r\times r}$ of $\rank (G(s)) \allowbreak =r$ such that
$N_1(s)G(s)=N_2(s)$. Moreover, $G(s)$ is polynomial by \cite[Main Theorem-4. p. 495]{Fo75}. From Theorem \ref{mincharact}, for every $s_0\in\overline{\FF}$,  $\rank(N_1(s_0)G(s_0))=\rank(N_2(s_0))=r$, which means that  $\rank(G(s_0))=r$.
By Remark \ref{remmincharact}, the matrix $G(s)$ is unimodular. The converse is immediate.
\hfill $\Box$

\hide{ {\bf Proof.} If
    $\mR(N_1(s))= \mR(N_2(s))$, there exists  $G(s)\in\FF(s)^{r\times r}$ of
$\rank (G(s))=r$
    such that
$N_1(s)G(s)=N_2(s)$.
We first see that $G(s)$ is polynomial.

Let $M(s)= \diag\left(\frac{\epsilon_1(s)}{\psi_1(s)},\dots,\frac{\epsilon_r(s)}{\psi_r(s)}\right)$  be the Smith--McMillan form of $G(s)$, i.e., there exists unimodular matrices $V(s)$ and $W(s)$ such that $G(s)=V(s)M(s)W(s)$. Denote $\hat N_1(s)=N_1(s)V(s)$ and $\hat N_2(s)=N_2(s)W(s)^{-1}$. Then, $\hat N_1(s)M(s)=\hat N_2(s)$. As $\hat N_2(s)$ is polynomial, the elements of the first column of $\hat N_1(s)$ are multiple of $\psi_1(s)$. Moreover,  if $\psi_1(s)\neq 1$, there exists $a\in \bar \FF$ such that $\psi_1(a)=0$, therefore, the first column of $\hat N_1(a)$ is zero, which is not possible because $\rank \hat N_1(a)=r$.

 Now, for every $s_0\in\overline{\FF}$,  $\rank(N_1(s_0)G(s_0))=\rank(N_2(s_0))=r$, which means that  $\rank(G(s_0))=r$.
    By Remark \ref{remmincharact}, the matrix $G(s)$ is unimodular.

    The converse is immediate.
\hfill $\Box$}

\medskip

The next proposition proves that there exist rational subspaces of $\FF (s)^p$ of dimension $r = 1, \ldots , p-1$ and arbitrary minimal indices. However, the $p$ minimal indices of $\FF (s)^p$ are all equal to zero.

\begin{proposition}\label{exminbasis}
  Given non negative integers $d_1, \dots, d_r$ and $q$,
there exists a minimal basis $M(s)\in  \FF[s]^{(r+q)\times r}$ with column degrees
  $d_1, \dots, d_r$ if and only if
$$
d_1=\dots=d_r=0 \text{ if }
q=0.
$$
\end{proposition}
{\bf Proof.}
Let $M(s)\in  \FF[s]^{(r+q)\times r}$ be a minimal basis with column degrees
$d_1, \dots, d_r$. If $q=0$, then the columns of $M(s)$ form a minimal basis of $\FF(s)^r$. The canonical basis, i.e., the columns of the identity matrix $I_r$, is another minimal basis of $\FF(s)^r$, because its order takes the least possible value equal to zero. Therefore, $d_1=\dots=d_r=0$.

\hide{$M(s)\in  \FF[s]^{r\times r}$ is unimodular and  column reduced with column degrees $d_1, \dots, d_r$.
As a consequence, $\deg(\det (M(s)))=0$ and $\deg(\det (M(s)))=\sum_{i=1}^rd_i$,  therefore, $d_1=\dots=d_r=0$.}

Conversely,
if $q=0$,  take $M(s)=I_r$; otherwise,
if $q>0$, take
$$
M(s)=\begin{bmatrix}
s^{ d_{1}}&0&0&\dots&0\\
1&s^{ d_{2}}&0&\dots&0\\
0&1&s^{d_{3}}&\dots&0\\
\vdots&\vdots&\vdots&\ddots&\vdots\\
0&0&0&\dots&s^{ d_r}\\
0&0&0&\dots&1\\
0&0&0&\dots&0\\
\vdots&\vdots&\vdots&\ddots&\vdots\\
0&0&0&\dots&0\\
\end{bmatrix}\in\FF[s]^{(r+q)\times r}.
$$
In both cases $M(s)\in  \FF[s]^{(r+q)\times r}$ is a minimal basis  with column degrees $d_1, \dots, d_r$ by Theorem \ref{mincharact}.  \hfill $\Box$

\begin{rem}
{\rm
Observe that if $q>0$, a minimal basis with  column degrees
$d_1, \dots, d_r$ always exists, for instance the matrix $M(s)$ in the proof. In fact, there are  other minimal bases with these column degrees because  if $P \in \FF^{(r+q) \times (r+q)}$ is any constant non singular matrix, then $P M(s)$ is a minimal basis by Theorem \ref{mincharact} and has the same column degrees as $M(s)$.

If $R(s) \in \FF(s)^{m\times n}$ has rank equal to $n$ ($m$), then its row-span (col-span) minimal indices are all equal to zero, because $\mR (R(s)^T) = \FF (s)^n$ ($\mR (R(s)) = \FF (s)^m$). Thus, these conditions must appear in all the results of this work about the prescription  of the row-span and col-span minimal indices of rational matrices.
\hfill $\Box$}
\end{rem}

Next, we summarize some properties of dual minimal bases and some consequences of these properties, which can be found in \cite{DeDoMaVa16, DmDoVa23, Fo75}.

\begin{definition}[{\rm \cite[Definition 2.10]{DeDoMaVa16}}]\label{defdual}
  Polynomial matrices $M(s)\in \FF[s]^{(r+q)\times r}$ and $N(s)\in
  \FF[s]^{(r+q)\times q}$ are said to be \textit{dual minimal bases} if they are minimal bases satisfying $M(s)^TN(s)=0$.
  \end{definition}

\begin{theorem}[{\rm \cite[Theorem 2.9]{DmDoVa23}}]\label{theodual}
Let $M(s)\in \FF[s]^{(r+q)\times r}$, $N(s)\in \FF[s]^{(r+q)\times q}$
be dual minimal bases with column degrees
$d_1, \dots, d_r$ and $d'_1, \dots, d'_q$, respectively. Then
\begin{equation}
  \label{eqsumdual}
  \sum_{i=1}^rd_i= \sum_{i=1}^qd'_i.
\end{equation}

Conversely, given  two lists of non negative integers $d_1, \dots, d_r$ and $d'_1, \dots, d'_q$,  satisfying \eqref{eqsumdual}, there
exists a pair of dual minimal bases
$M(s)\in \FF[s]^{(r+q)\times r}$ and $N(s)\in \FF[s]^{(r+q)\times q}$
with precisely these column
degrees, respectively.
\end{theorem}

The next result is an immediate corollary of Theorem \ref{theodual} and has been stated and proved in \cite[Corollary 2.10]{DmDoVa23} for polynomial matrices. The proof in \cite{DmDoVa23} remains valid for rational matrices.
\begin{corollary} \label{cormin}
  Let $A(s)\in \FF(s)^{m\times n}$, $\rank(A(s))=r$, with col-span minimal indices $k_1,\dots, k_r$,  row-span minimal indices $\ell_1,\dots,\ell_r$, right minimal indices $d_1, \dots, d_{n-r}$,
 and left minimal indices $v_1, \dots, v_{m-r}$.
Then
\begin{equation}\label{eqsums}
\sum_{i=1}^{m-r}v_i= \sum_{i=1}^rk_i   \quad \mbox{ and } \quad
\sum_{i=1}^{n-r}d_i= \sum_{i=1}^r\ell_i.
 \end{equation}
\end{corollary}

The next result is specific for polynomial matrices. It is an immediate consequence of Theorem  \ref{theoexistenceDeDoVa152} and Corollary \ref{cormin}.

\begin{corollary}\label{corxy}
    Let $A(s)\in\FF[s]^{m\times n}$, $\rank(A(s))=r$, $\deg(A(s))=d$, with col-span minimal indices $k_1,\dots, k_r$, row-span minimal indices $\ell_1, \dots,\ell_r$, invariant factors
    $\alpha_1(s)\mid \dots \mid \alpha_r(s)$, and  $f_1, \dots, f_r$ as partial multiplicities of $\infty$.
    Then,
    $f_1=0$ and
\begin{equation}\label{eqsumklfa}
\sum_{i=1}^{r}k_i+\sum_{i=1}^{r}\ell_i+\sum_{i=1}^{r}f_i+\sum_{i=1}^{r}\deg(\alpha_i)=rd.
  \end{equation}
\end{corollary}

\section{Statement of the problems for polynomial matrices} \label{problems}

In this section, we present a formal statement of the problems studied in this paper for polynomial matrices. They are solved in Section \ref{main}.  For brevity, we avoid the statement of the corresponding problems for rational matrices, which are solved in Section \ref{sec.rational} as a consequence of the results obtained for polynomial matrices.

\begin{problem}\label{probkl2}
 Let $m$, $n$, $r\leq \min\{m,n \}$ be  positive integers and $d$ a non negative integer. Let $\alpha_1(s)\mid \dots \mid \alpha_r(s)$ be monic polynomials in $\FF[s]$,
$(f_r, \ldots, f_{1})$   a  partition of non negative integers, and let $K(s)\in \FF[s]^{m\times r}$,  $L(s)^T\in \FF[s]^{n\times r}$ be minimal bases.
Find necessary and sufficient conditions for the existence of a polynomial matrix $A(s)\in\FF[s]^{m\times n}$ of   $\rank(A(s))= r$,
$\deg(A(s))= d$, with
   $\alpha_1(s),\dots, \alpha_r(s)$ as invariant factors,  $f_1, \dots, f_r$ as partial multiplicities of $\infty$, $\mR(A(s))=\mR(K(s))$ and $\mR(A(s)^T)=\mR(L(s)^T)$.
\end{problem}

\begin{problem}\label{probkl}
 Let $m$, $n$, $r\leq \min\{m,n \}$ be  positive integers and $d$  a non negative integer.
Let $\alpha_1(s)\mid \dots \mid \alpha_r(s)$ be monic polynomials in $\FF[s]$ and
let $(f_r, \ldots, f_{1})$,  $(k_1, \ldots, k_{r})$, $(\ell_1, \ldots, \ell_{r})$ be partitions of non negative integers.
Find necessary and sufficient conditions for the existence of a polynomial matrix $A(s)\in\FF[s]^{m\times n}$ of $\rank(A(s))= r$, $\deg(A(s))= d$, with $\alpha_1(s),\dots, \alpha_r(s)$ as invariant factors, $f_1, \dots, f_r$ as partial multiplicities of $\infty$, $k_1,\dots, k_r$ as col-span minimal indices, and $\ell_1,\dots,\ell_r$ as row-span minimal indices.
\end{problem}
\begin{problem}\label{probklcv}
 Let $m$, $n$, $r\leq \min\{m,n \}$ be  positive integers and $d$  a non negative integer.
Let $\alpha_1(s)\mid \dots \mid \alpha_r(s)$ be monic polynomials in $\FF[s]$ and
let $(f_r, \ldots, f_{1})$, $(k_1, \ldots, k_{r})$, $(\ell_1, \ldots, \ell_{r})$, $(d_1, \dots, d_{n-r})$,  $(v_1, \dots, v_{m-r})$  be partitions of non negative integers.
Find necessary and sufficient conditions for the existence of a polynomial matrix $A(s)\in\FF[s]^{m\times n}$ of  $\rank(A(s))= r$, $\deg(A(s))= d$, with $\alpha_1(s),\dots, \alpha_r(s)$ as invariant factors, $f_1, \dots, f_r$ as partial multiplicities of $\infty$, $k_1,\dots, k_r$ as col-span minimal indices, $\ell_1,\dots,\ell_r$ as row-span minimal indices, and right and left minimal indices $d_1, \dots, d_{n-r}$ and $v_1, \dots, v_{m-r}$, respectively.
\end{problem}

Problems \ref{probkl2}, \ref{probkl},  and \ref{probklcv} are solved in Theorems \ref{theomainlk2}, \ref{cormainlk}, and \ref{cormainlkcv}, respectively. We will see that the solution of Problem \ref{probkl} follows immediately from combining that of Problem \ref{probkl2} with Proposition \ref{exminbasis}. The solution of Problem \ref{probklcv} requires to combine the solutions of Problems \ref{probkl2} and \ref{probkl} with the non-trivial results stated in Theorem \ref{theodual} and Corollary \ref{cormin}. Thus, the key original contribution of this manuscript is the solution of Problem \ref{probkl2}.

The counterparts of Problems \ref{probkl2}, \ref{probkl},  and \ref{probklcv} for rational matrices are solved in Theorems \ref{theomainlk2_rat}, \ref{cormainlk_rat}, and \ref{cormainlkcv_rat}, respectively.

\section{Factorizations of polynomial matrices and a reformulation of Problem \ref{probkl2}} \label{sec.factorization}
The goal of this section is to establish Theorem \ref{theoreformulation}, which plays a key role in the solution of Problem \ref{probkl2}. Theorem \ref{theoreformulation} shows that Problem \ref{probkl2} is equivalent to an existence problem involving only the invariant factors of a regular polynomial matrix and the partial multiplicities of $\infty$  of a two-sided diagonal scaling of it. The proof of Theorem \ref{theoreformulation} requires the development of several results on factorizations of polynomial matrices. The first one of them proves that column proper polynomial matrices with full column rank admit the factorization given in Theorem \ref{thm.factcolproper}.

\begin{theorem} \label{thm.factcolproper} Let $K(s)\in\FF[s]^{m\times r}$, $\rank(K(s))=r$,  be column proper with column degrees $k'_1,\dots, k'_r$. Then, there exists a biproper matrix $B(s) \in\FF_{pr}(s)^{m\times m}$ such that
 \begin{equation*}\label{eq.factcolproper}
   K(s) = B(s)
   \left[\begin{array}{c}\diag(s^{k'_1},\ldots,s^{k'_r})\\0\end{array}\right].
 \end{equation*}
\end{theorem}

\textbf{Proof.} Let us express $K(s)$ as
$$
K(s) = K_h \diag(s^{k'_1},\ldots,s^{k'_r}) + L(s),
$$
where $K_h \in \FF^{m\times r}$ is the highest column degree coefficient matrix of $K(s)$ and has rank equal to $r$, and the degree of the $i$-th column of $L(s)  \in\FF[s]^{m\times r}$ is strictly less than $k'_i$ for $i = 1, \ldots , r$. Therefore,
$$
K(s) = (K_h + \widehat{L}(s)) \, \diag(s^{k'_1},\ldots,s^{k'_r}),
$$
with $\widehat{L}(s) = L(s) \, (\diag(s^{k'_1},\ldots,s^{k'_r}))^{-1}$ strictly proper. Let $\widetilde{K}_h \in \FF^{m\times (m-r)}$ be any constant matrix such that $\left[\begin{array}{cc} K_h & \widetilde{K}_h \end{array}\right]\in \FF^{m\times m}$ is invertible. Note that
$$
   K(s) = \left[\begin{array}{cc} K_h + \widehat{L}(s) & \widetilde{K}_h \end{array}\right]
   \left[\begin{array}{c}\diag(s^{k'_1},\ldots,s^{k'_r})\\0\end{array}\right].
$$
We prove that $B(s) = \left[\begin{array}{cc} K_h + \widehat{L}(s) & \widetilde{K}_h \end{array}\right]$ is biproper by proving that $\det (B(s))$ is a non zero biproper rational scalar function.
 Since the determinant of a matrix is a multilinear function of its columns, we can express $\det (B(s))$ as
$$
\det (B(s)) = \det \left(\left[\begin{array}{cc} K_h & \widetilde{K}_h \end{array}\right]\right) + g(s),
$$
where $g(s)$ is a sum of $2^r - 1$ determinants each containing at least one column of $\widehat{L}(s)$. Therefore, each one of these determinants is a strictly proper rational function (use a cofactor expansion through a column of $\widehat{L}(s)$). Hence, $\det (B(s)) = c + g(s)$  where $0 \ne c \in \FF$ and $g(s)$ is strictly proper. This implies that $\det (B(s))$ is biproper.
\hfill $\Box$

\medskip

The second result of this section deals with polynomial matrices which are products of three factors, the first and the third ones satisfying certain properties which, in particular, are fulfilled by minimal bases.

\begin{lemma}\label{lemm.sinvfactcolprop}
  Let $K(s)\in \FF[s]^{m\times r }$ and  $L(s)\in \FF[s]^{r\times n }$ be polynomial matrices  with $\rank(K(s))=\rank(L(s))=r$,  and let $k'_1,\dots, k'_r$ be the column degrees of $K(s)$  and $\ell'_1,\dots, \ell'_r$ be the row degrees of $L(s)$. Let $\widehat E(s)\in \FF[s]^{r\times r }$ and
  $A(s)=K(s)\widehat E(s)L(s)$.
\begin{enumerate}
\item[\rm (i)] If the invariant factors of $K(s)$ and $L(s)$ are  equal to $1$, then
       $A(s)$ and $\widehat E(s)$ have the same invariant factors.

\item[\rm (ii)] If $K(s)$ and $L(s)^T$ are column proper, then $$
F(s) = \diag(s^{k'_1}, \dots, s^{k'_r})  \widehat E(s)\diag(s^{\ell'_1}, \dots, s^{\ell'_r})
$$
has the same degree and partial multiplicities of $\infty$ as $A(s)$.
\end{enumerate}
\end{lemma}

\textbf{Proof.} To prove part (i) note that, by Remark \ref{remmincharact}, $\rank K(s_0)=r$ and $\rank L(s_0)=r$ for all $s_0\in\overline{\efe}$.
Then, by \cite[Lemma 2.16-(b)]{DeDoVa15} there exist $K'(s)\in \FF[s]^{m\times (m-r) }$ and  $L'(s)\in \FF[s]^{(n-r)\times n }$ such that
$U(s)=\begin{bmatrix}K(s)&K'(s)\end{bmatrix}\in \FF[s]^{m\times m }$ and $V(s)=\begin{bmatrix}L(s)\\L'(s)\end{bmatrix}\in \FF[s]^{n\times n }$
 are unimodular.
Let
$\widehat E'(s)=\begin{bmatrix}\widehat E(s)&0\\0&0\end{bmatrix}\in  \FF[s]^{m\times n }$. Then $\widehat E'(s)$ has the same invariant factors as $\widehat E(s)$ and
$
U(s)\widehat E'(s)V(s)=A(s)
$;
hence $A(s)$ and $\widehat E(s)$ have the same invariant factors.

To prove part (ii), we use Theorem \ref{thm.factcolproper} to express
$$
K(s) = B_1(s)
   \left[\begin{array}{cc}\diag(s^{k'_1},\ldots,s^{k'_r})\\0\end{array}\right],
\;
L(s) =
\left[\begin{array}{cc}\diag(s^{\ell'_1},\ldots,s^{\ell'_r})&0\end{array}\right] B_2(s),$$
where $B_1(s)\in\FF_{pr}(s)^{m\times m}, B_2(s)\in\FF_{pr}(s)^{n\times n}$ are  biproper matrices. Then,
\begin{align*}
  A(s) & =  B_1(s) \left[\begin{array}{cc}\diag(s^{k'_1},\ldots,s^{k'_r})\\0\end{array}\right] \widehat{E} (s) \left[\begin{array}{cc}\diag(s^{\ell'_1},\ldots,s^{\ell'_r})&0\end{array}\right]
   B_2(s)
    \\
   & = B_1(s) \left[\begin{array}{cc}F(s)&0\\0&0\end{array}\right] B_2 (s) .
\end{align*}
Thus, $A(s)$ and $\left[\begin{array}{cc}F(s)&0\\0&0\end{array}\right]$ are equivalent at  infinity; hence they have the same degree and the same partial
 multiplicities of $\infty$ (see  the discussion after equation \eqref{eqqf}).
\hfill $\Box$

\medskip

\textit{Minimal rank factorizations} of polynomial matrices were introduced in \cite[Definition 3.12, Theorem 3.11]{DmDoVa23}. Some properties of minimal rank factorizations into three factors are presented in the next theorem. Observe that Theorem \ref{theofactDmDoDoEXPAND}-(iii) proves  a new property.

\begin{theorem} \label{theofactDmDoDoEXPAND}
Every  polynomial matrix $A(s)\in \FF[s]^{m\times n}$ of rank $r$ can be factorized as
$$A(s)=K(s)\widehat E(s)L(s), \quad K(s)\in \FF[s]^{m\times r}, \quad \widehat E(s)\in \FF[s]^{r\times r}, \quad
L(s)\in \FF[s]^{r\times n},$$
where
\begin{enumerate}
\item[\rm (i)] the columns of $K(s)$  form a minimal basis of $\mR(A(s))$,
the columns of $L(s)^T$  form a minimal basis  of $\mR(A(s)^T)$,
\item[\rm (ii)] $\widehat E(s)$ has the same  invariant factors as $A(s)$, and
\item[\rm (iii)] if the column degrees of $K(s)$ are $k'_1,\dots, k'_r$ and
the column degrees of $L(s)^T$ are $\ell'_1,\dots, \ell'_r$, then
$$
F(s) = \diag(s^{k'_1}, \dots, s^{k'_r}) \, \widehat E(s)\, \diag(s^{\ell'_1}, \dots, s^{\ell'_r})
$$
has the same degree and the same partial multiplicities of $\infty$ as $A(s)$.
\end{enumerate}
\end{theorem}

\textbf{Proof.}
 The existence of the factorization and parts (i) and (ii) are in \cite[Theorem 3.11-(i)]{DmDoVa23}. Part (iii) follows from Lemma \ref{lemm.sinvfactcolprop}-(ii).
\hfill $\Box$

  \begin{rem}\label{remfactDmDoDo}\, \textit{ (On the minimal bases of Theorem \ref{theofactDmDoDoEXPAND})}

{\rm
The minimal bases of $\mR(A(s))$ and $\mR(A(s)^T)$ can be arbitrarily chosen in Theorem \ref{theofactDmDoDoEXPAND}. Let the columns of $K'(s)\in \FF[s]^{m\times r}$ and  $L'(s)^T\in\FF[s]^{n\times r}$ form minimal bases of  $\mR(A(s))$ and $\mR(A(s)^T)$, respectively. By Proposition \ref{eqbasis}, there exist unimodular matrices $U(s), V(s)\in \FF[s]^{r\times r}$ such that $K(s)=K'(s)U(s)$ and $L(s)=V(s)L'(s)$. Then $$A(s)=K(s)\widehat E(s)L(s)=K'(s)\widehat E'(s)L'(s),$$
where $\widehat E'(s)=U(s)\widehat E(s)V(s)$ has the same invariant factors as $\widehat E(s)$,
and the matrix $F'(s)$  defined from the factorization $A(s) = K'(s)\widehat E'(s)L'(s)$  as in Theorem \ref{theofactDmDoDoEXPAND}-(iii) has the same degree and the same partial multiplicities of $\infty$ as $A(s)$ by Lemma \ref{lemm.sinvfactcolprop}-(ii).

In particular, the columns of $K(s)$ and the rows of $L(s)$
can be arbitrarily ordered. \hfill $\Box$}
\end{rem}

Next, we prove the announced reformulation of Problem \ref{probkl2}.

\begin{theorem}\label{theoreformulation}
Let $m$, $n$, $r\leq \min\{m,n \}$ be positive integers and $d$ a non negative integer.
Let $\alpha_1(s)\mid \dots \mid \alpha_r(s)$ be monic polynomials in $\FF [s]$,  let
 $(f_r, \ldots, f_{1})$ be a  partition of non negative integers, and let $K(s)\in \FF[s]^{m\times r}$,
$L(s)^T\in \FF[s]^{n\times r}$ be minimal bases with column degrees $k_1\geq \dots \geq  k_r$ and $\ell_1\geq \dots\geq  \ell_r$, respectively.
The following statements are equivalent:

\begin{itemize}
  \item[\rm (i)]
There exists a polynomial matrix $A(s)\in\FF[s]^{m\times n}$ with  $\rank(A(s))= r$,
$\deg(A(s))= d$, $\alpha_1(s), \dots, \alpha_r(s)$ as invariant factors, $f_1, \dots, f_r$ as partial multiplicities of $\infty$,
$\mR(A(s))=\mR(K(s))$, and $\mR(A(s)^T)=\mR(L(s)^T)$.
\item[\rm (ii)]
There exists $  E(s)\in \FF[s]^{r\times r}$ with   $\alpha_1(s), \dots, \alpha_r(s)$ as invariant factors,
such that the polynomial matrix
 \begin{equation}\label{eqmatF}
F(s)=\diag(s^{k_1}, \dots, s^{k_r})  E(s)\diag(s^{\ell_1}, \dots, s^{\ell_r})
  \end{equation}
has $f_1, \dots, f_r$ as partial multiplicities of $\infty$ and  $\deg(F(s))=d$.
  \end{itemize}
  \end{theorem}

{\bf Proof.}\,

\underline{(i) $\Rightarrow$ (ii).}
It follows from  Theorem \ref{theofactDmDoDoEXPAND} and Remark \ref{remfactDmDoDo}.

\underline{(ii) $\Rightarrow$ (i)}.
Let us assume that (ii) holds. Let $A(s)=K(s)E(s)L(s) \in \FF[s]^{m \times n}$. By Remark \ref{remmincharact} and Lemma \ref{lemm.sinvfactcolprop}-(i),  the invariant factors of $A(s)$ are  $\alpha_1(s), \dots, \alpha_r(s)$. This implies that $\rank (A(s)) =r$. Moreover,  we have $\mR(A(s)) \subseteq \mR(K(s))$ and $\mR(A(s)^T) \subseteq \mR(L(s)^T)$. Since each of these four subspaces has dimension $r$, we get  $\mR(A(s)) = \mR(K(s))$ and $\mR(A(s)^T) = \mR(L(s)^T)$.
 By Lemma \ref{lemm.sinvfactcolprop}-(ii), $\deg(A(s))=d$ and $A(s)$ has $f_1, \dots, f_r$ as  partial multiplicities of $\infty$.
   \hfill $\Box$

\begin{rem}\label{remreformulationord}
{\rm
In item (ii) of Theorem \ref{theoreformulation}, the indices $k_1, \ldots, k_r$ and $\ell_1, \ldots, \ell_r$ can be arbitrarily ordered. That is to say, if
  $k'_1, \dots, k'_r$ and $\ell'_1, \dots, \ell'_r$ are reorderings of
  $k_1, \dots, k_r$ and $\ell_1, \dots, \ell_r$, respectively,
  then
item (ii) of Theorem \ref{theoreformulation} holds if and only if
there exists $  E'(s)\in \FF[s]^{r\times r}$ with  invariant factors
$\alpha_1(s), \dots, \alpha_r(s)$, such that the matrix
 \begin{equation*}\label{eqmatFp}
F'(s)=\diag(s^{k'_1}, \dots, s^{k'_r})  E'(s)\diag(s^{\ell'_1}, \dots, s^{\ell'_r})
  \end{equation*}
  has $f_1, \dots, f_r$ as partial multiplicities of $\infty$ and  $\deg(F'(s))=d$.

 To prove it, assume that item (ii) of Theorem \ref{theoreformulation} holds and take
 $P, Q\in \FF^{r\times r}$  permutation matrices such that
 $$
 P\diag(s^{k_1}, \dots, s^{k_r})P^T=\diag(s^{k'_1}, \dots, s^{k'_r}),$$$$
Q\diag(s^{\ell_1}, \dots, s^{\ell_r})Q^T=\diag(s^{\ell'_1}, \dots, s^{\ell'_r}),
 $$
and let $E'(s)= PE(s)Q^T$. Then the invariant factors of $E'(s)$ are
$\alpha_1(s),\dots,\alpha_r(s)$ and
$$
PF(s)Q^T
=\diag(s^{k'_1}, \dots, s^{k'_r})  E'(s)\diag(s^{\ell'_1}, \dots, s^{\ell'_r})=F'(s).
$$
Therefore, $F(s)$ and $F'(s)$ are equivalent at $\infty$; hence they have the same degree and the same partial multiplicities of $\infty$. The proof of the converse is analogous.}
  \end{rem}

\section{Solution to the problems for polynomial matrices}
\label{main}
In this section we provide necessary and sufficient conditions for the existence problems posed in Section \ref{problems}. We emphasize that the conditions found are sufficient only if the underlying field is algebraically closed; Example \ref{exnoac} illustrates that this constraint is essential. In contrast, they are necessary over arbitrary fields.  Theorem \ref{theomainlk2} solves Problem \ref{probkl2}.

\begin{theorem}\label{theomainlk2}
Let $m$, $n$, $r\leq \min\{m,n \}$ be  positive integers and $d$ a non negative integer.
Let
$\alpha_1(s)\mid \dots \mid \alpha_r(s)$ be monic polynomials in $\FF[s]$,
$(f_r, \ldots, f_{1})$  a  partition of non negative integers and  $K(s)\in \FF[s]^{m\times r}$,  $L(s)^T\in \FF[s]^{n\times r}$  minimal  bases,
with column degrees $k_1\geq \dots\geq  k_r$ and $\ell_1\geq \dots\geq  \ell_r$, respectively.
Let
$g_1\geq \dots \geq g_r$ be  the decreasing reordering of $k_r+\ell_1, \dots, k_1+\ell_r$.

If $A(s)\in\FF[s]^{m\times n}$ is a polynomial matrix with $\rank(A(s))=r$, $\deg(A(s))=d$,
$\alpha_1(s),\dots, \alpha_r(s)$ as invariant factors,  $f_1, \dots, f_r$ as partial multiplicities of $\infty$, $\mR(A(s))  = \mR(K(s))$ and $\mR(A(s)^T)=\mR(L(s)^T)$, then
\begin{equation}\label{eqf1}
    f_1=0
\end{equation}
and
\begin{equation}\label{eqprec}
  (d-g_r, \dots, d-g_1)\prec (\deg(\alpha_r)+f_r, \dots, \deg(\alpha_1)+f_1).
  \end{equation}

Conversely, if $\FF$ is  algebraically closed, conditions \eqref{eqf1} and \eqref{eqprec}
are sufficient for the existence of a polynomial matrix $A(s)\in\FF[s]^{m\times n}$
 with $\rank(A(s))=r$, $\deg(A(s))=d$,
$\alpha_1(s),\dots, \alpha_r(s)$ as invariant factors,  $f_1, \dots, f_r$ as partial multiplicities of $\infty$, $\mR(A(s)) \allowbreak =\mR(K(s))$ and $\mR(A(s)^T)=\mR(L(s)^T)$.
\end{theorem}

Despite the simplicity of the statement, the proof of Theorem \ref{theomainlk2} is long and involved, and is postponed. The proof of the necessity of this theorem is given in Subsection \ref{subsecnec}, while the sufficiency is proved in Subsection \ref{subsecsuf}. The proof of the necessity uses a determinantal lemma, Lemma \ref{lemmamingID}, which is presented together with its long proof in Appendix \ref{secappend}.

\begin{rem}\label{remprec} \textit{ (On the conditions \eqref{eqf1} and \eqref{eqprec} in Theorem \ref{theomainlk2}).}
{\rm
\begin{enumerate}
\item The strong, non trivial condition in Theorem \ref{theomainlk2} is \eqref{eqprec}, because condition \eqref{eqf1} was expected (recall the discussion after \eqref{eqqf}).

\item Condition \eqref{eqprec} amounts to $r$ conditions due to the definition of majorization. The last of such conditions, i.e., the equality of the sums of the terms of the two sequences involved in \eqref{eqprec}, is precisely \eqref{eqsumklfa}, which in turn is equivalent to the condition \eqref{eqIST}, i.e., to the index sum theorem constraint.

\item Using the observation above, condition \eqref{eqprec} is equivalent to conditions \eqref{eqsumklfa} and \begin{equation*}\label{eqcondleq}
\sum_{i=1}^k(g_i+\deg(\alpha_i)+f_i)\leq kd, \quad 1\leq k\leq r-1.
\end{equation*}

\item For $k=1$, the inequality above yields $\deg(\alpha_1)+f_1 \leq d - g_1$. This implies that if \eqref{eqprec} holds, then all the terms of the sequence in the left hand side of \eqref{eqprec} are non negative.

\item Theorem \ref{theomainlk2} can be used to determine the existence of a polynomial matrix without eigenvalues (neither finite nor infinite) and prescribed degree, rank and column and row spaces. For this purpose, we must take $\alpha_1(s)=\dots= \alpha_r(s)=1$ and  $f_1= \dots= f_r=0$. In this case, \eqref{eqprec} reduces to
  $$
\ell_i+k_{r-i+1}=d, \quad 1\leq i \leq r.
  $$
This condition is related to \cite[Theorem 3.19]{DmDoVa23}. \hfill $\Box$
\end{enumerate}}
\end{rem}

In Theorems \ref{cormainlk} and \ref{cormainlkcv} we give solutions to Problems \ref{probkl} and  \ref{probklcv}, respectively. The proofs of both theorems are based on Theorem \ref{theomainlk2}. Theorem \ref{cormainlk} is the counterpart of Theorem \ref{theoexistenceDeDoVa152} when the col-span and row-span minimal indices are prescribed instead of the right and left minimal indices. Observe that this change in the prescribed data leads to a considerable increase in the complexity of the required conditions. The key condition in Theorem \ref{cormainlk} is \eqref{eqprec}, while \eqref{eqx>0} and \eqref{eqy>0} simply take into account that if the rank of an $m \times n$ polynomial matrix $A(s)$ is $n$ ($m$), then $\mathcal{R} (A(s)^T)= \FF(s)^{n}$ ($\mathcal{R} (A(s))= \FF(s)^{m}$) and the minimal indices of $\FF(s)^{n}$ ($\FF(s)^{m}$) are zero.

\begin{theorem}\label{cormainlk}
  Let $m$, $n$, $r\leq \min\{m,n \}$ be  positive integers and $d$  a non negative integer.
Let $\alpha_1(s)\mid \dots \mid \alpha_r(s)$ be monic polynomials  in $\FF [s]$ and
$(f_r, \ldots, f_{1})$,  $(k_1, \ldots, k_{r})$, $(\ell_1, \ldots, \ell_{r})$ partitions of non negative integers.
Let
$g_1\geq \dots \geq g_r$ be  the decreasing reordering of $k_r+\ell_1, \dots, k_1+\ell_r$.

If $A(s)\in\FF[s]^{m\times n}$ is a polynomial matrix with $\rank(A(s))=r$, $\deg(A(s))=d$,
$\alpha_1(s),\dots, \alpha_r(s)$ as invariant factors,  $f_1, \dots, f_r$ as partial multiplicities of $\infty$, $k_1,\dots, k_r$ as col-span minimal indices, and $\ell_1, \dots,\ell_r$ as row-span minimal indices, then \eqref{eqf1}, \eqref{eqprec},
\begin{equation}\label{eqx>0}
\ell_1 = \cdots = \ell_r =0 \text{ if } r = n,
\end{equation}
and
\begin{equation}\label{eqy>0}
k_1 = \cdots = k_r =0 \text{ if } r = m.
\end{equation}

Conversely, if $\FF$ is  algebraically closed,  conditions \eqref{eqf1}, \eqref{eqprec}, \eqref{eqx>0} and \eqref{eqy>0}
are sufficient for the existence of a polynomial matrix $A(s)\in\FF[s]^{m\times n}$
 with $\rank(A(s))=r$, $\deg(A(s))=d$,
$\alpha_1(s),\dots, \alpha_r(s)$ as invariant factors,  $f_1, \dots, f_r$ as partial multiplicities of $\infty$, $k_1,\dots, k_r$ as col-span minimal indices, and $\ell_1, \dots,\ell_r$ as row-span minimal indices.
\end{theorem}

{\bf Proof.}
By Proposition \ref{exminbasis},  there exist minimal bases
$K(s)\in \FF[s]^{m\times r}$,
$L(s)^T\in \FF[s]^{n\times r}$, with column degrees
$k_1, \dots,  k_r$, and $\ell_1, \dots,  \ell_r$, respectively, if and only if  (\ref{eqx>0}) and
(\ref{eqy>0}) hold. Then, the result follows from Theorem \ref{theomainlk2}.
\hfill $\Box$

\medskip

If $\FF$ is not  algebraically closed, conditions \eqref{eqf1}, \eqref{eqprec}, \eqref{eqx>0} and \eqref{eqy>0} of Theorem \ref{cormainlk}
are not sufficient, as we can see in the following example. The same happens with the conditions in Theorems \ref{theomainlk2} and \ref{cormainlkcv}.

\begin{example}
\label{exnoac}
{\rm
Let $\FF=\RR$, $m=n=3$, $r=2$, $d=7$, $\alpha_1(s)=1$, $\alpha_2(s)=(s^2+1)^2$, $f_1=f_2=0$, $k_1=5$, $k_2=0$, $\ell_1=4$, $\ell_2=1$.
Then $r<m$, $r<n$, $(g_1, g_2)= (k_1+\ell_2, k_2+\ell_1)=(6,4)$ and $(d-g_2, d-g_1)=(3,1)\prec(4, 0)=(\deg(\alpha_2)+f_2, \deg(\alpha_1)+f_1)$. Thus, \eqref{eqf1}, \eqref{eqprec}, \eqref{eqx>0} and \eqref{eqy>0} hold.
If there exists $A(s)\in \RR[s]^{3\times 3}$, with $\rank (A(s)) = 2$, $\deg(A(s))=7$,
$\alpha_1(s)=1$, $\alpha_2(s)=(s^2+1)^2$ as invariant factors, $f_1=f_2=0$ as partial multiplicities of $\infty$, $5,0$ as col-span minimal indices, and $4,1$ as row-span minimal indices, then, by Theorem \ref{theoreformulation}, there exists a polynomial matrix
$E(s)=\begin{bmatrix}e_{1,1}(s)&e_{1,2}(s)\\e_{2,1}(s)&e_{2,2}(s)\end{bmatrix}\in \RR[s]^{2\times 2}$ with invariant factors $\alpha_1(s)=1, \alpha_2(s)=(s^2+1)^2$
and such that
$$
F(s)=\diag(s^{5},s^{0})  E(s)\diag(s^{4},s^{1})=\begin{bmatrix}s^{9} e_{1,1}(s)&s^{6} e_{1,2}(s)\\s^{4} e_{2,1}(s)&s e_{2,2}(s)\end{bmatrix}
$$
has degree  $7$ and $f_1=f_2=0$ as partial multiplicities of $\infty$.  Then,
$$
e_{1,1}(s)=0, \quad \deg(e_{1,2} (s))\leq 1, \quad \deg(e_{2,1} (s))\leq 3, \quad \deg(e_{2,2} (s))\leq 6.
$$
Consequently,
$\alpha_2(s)=(s^2+1)^2= c \, e_{1,2}(s) \, e_{2,1}(s)$,  where  $\deg(e_{1,2} (s))=1$, $\deg(e_{2,1} (s))=3$, and $c\ne 0$ is a real number. This is a contradiction because $(s^2+1)^2$ has no divisors of degree $1$ over $\mathbb{R}[s]$.
}
\end{example}

Our last result in this section, Theorem \ref{cormainlkcv}, combines all the data that are prescribed in Theorems \ref{theoexistenceDeDoVa152} and \ref{cormainlk}, since, in addition to the invariant factors and the partial multiplicities of $\infty$, the four sequences of minimal indices are prescribed. Observe that Theorem \ref{cormainlkcv} requires the extra condition \eqref{eqsums} with respect to Theorem \ref{cormainlk} and that \eqref{eqsums} implies \eqref{eqx>0} and \eqref{eqy>0}.

\begin{theorem}\label{cormainlkcv}
Let $m$, $n$, $r\leq \min\{m,n \}$ be  positive integers and $d$  a non negative integer.
Let $\alpha_1(s)\mid \dots \mid \alpha_r(s)$ be monic polynomials  in $\FF[s]$ and
$(f_r, \ldots, f_{1})$, $(k_1, \ldots, k_{r})$, $(\ell_1, \ldots, \ell_{r})$, $(d_1, \dots, d_{n-r})$,  $(v_1, \dots, v_{m-r})$  partitions of non negative integers.
Let
$g_1\geq \dots \geq g_r$ be  the decreasing  reordering of $k_r+\ell_1, \dots, k_1+\ell_r$.

If $A(s)\in\FF[s]^{m\times n}$ is a polynomial matrix with $\rank(A(s))=r$, $\deg(A(s))=d$,
$\alpha_1(s),\dots, \alpha_r(s)$ as invariant factors,  $f_1, \dots, f_r$ as partial multiplicities of $\infty$, $k_1,\dots, k_r$ as col-span minimal indices, $\ell_1, \dots,\ell_r$ as row-span minimal indices, and
$d_1, \dots, d_{n-r}$ and $v_1, \dots, v_{m-r}$ as right and left minimal indices, respectively, then \eqref{eqsums}, \eqref{eqf1} and \eqref{eqprec}  hold.

Conversely, if $\FF$ is  algebraically closed, conditions \eqref{eqsums},
\eqref{eqf1} and \eqref{eqprec} are sufficient for the existence of a polynomial matrix  $A(s)\in\FF[s]^{m\times n}$, with $\rank(A(s))=r$, $\deg(A(s))=d$,  $\alpha_1(s), \dots, \alpha_r(s)$ as invariant factors,
$f_1, \dots, f_r$ as partial multiplicities of $\infty$, $k_1,\dots, k_r$ as col-span minimal indices,  $\ell_1, \dots,\ell_r$ as row-span minimal indices,  and
$d_1, \dots, d_{n-r}$ and  $v_1, \dots, v_{m-r}$ as right and left minimal indices, respectively.
\end{theorem}
{\bf Proof.}
The necessity follows from Corollary \ref{cormin} and Theorem \ref{cormainlk}.

Assume now that $\FF$ is  algebraically closed and conditions \eqref{eqsums},
\eqref{eqf1} and \eqref{eqprec} hold. If $m=r$, then from  \eqref{eqsums} we obtain $\sum_{i=1}^rk_i=\sum_{i=1}^0v_i=0$, i.e., $k_1=\dots= k_r=0$, and by Proposition \ref{exminbasis}  there exists a minimal basis $K(s)\in \FF[s]^{m\times r}$ with column degrees $k_1,\dots, k_r$. If $m>r$, from (\ref{eqsums}), by Theorem \ref{theodual} there exist minimal bases $K(s)\in \FF[s]^{m\times r}$, $N_K (s)\in \FF[s]^{m\times (m-r)}$ with column degrees $k_1,\dots, k_r$ and $v_1,\dots, v_{m-r}$, respectively, such that $N_K(s)^TK(s)=0$.
Analogously, using \eqref{eqsums}, Proposition \ref{exminbasis}  and Theorem \ref{theodual}, there exists a minimal basis $L(s)^T\in \FF[s]^{n\times r}$ with column degrees $\ell_1,\dots, \ell_r$, and if $n>r$ there exists a minimal basis $N_L(s)\in \FF[s]^{n\times (n-r)}$
with  column degrees $d_1, \dots d_{n-r}$ such that  $L(s)^TN_L(s)=0$.

By Theorem \ref{theomainlk2}, from (\ref{eqf1}) and (\ref{eqprec}),
there exists a  polynomial matrix $A(s)\in\FF[s]^{m\times n}$, with $\rank(A(s))=r$, $\deg(A(s))=d$, $\alpha_1(s),\dots,\alpha_r(s)$ as invariant factors,  $f_1,\dots, f_r$ as partial multiplicities of $\infty$, $\mR(A(s))=\mR(K(s))$ and $\mR(A(s)^T)=\mR(L(s)^T)$, and therefore with $k_1, \ldots , k_r$ as col-span minimal indices and $\ell_1, \ldots , \ell_r$ as row-span minimal indices.

If $m>r$, since $N_K(s)^TK(s)=0$, the matrix $N_K(s)$ is a minimal basis of
$\mathcal{N}_\ell(A(s))$, therefore the left  minimal indices of $A(s)$ are  $v_1, \dots, v_{m-r}$. Analogously, if $n>r$, $d_1, \dots, d_{n-r}$ are the right minimal indices of $A(s)$.
\hfill $\Box$

\subsection{Proof of Theorem \ref{theomainlk2}: necessity  of the conditions}
\label{subsecnec}

In this subsection we prove the necessity part of Theorem \ref{theomainlk2}.
The proof  requires to select some minors of a polynomial matrix satisfying certain properties. To do it, we need to introduce some notation and  a technical lemma.
In order  not to pause the proof of the theorem,  the proof of the technical lemma is postponed to Appendix \ref{secappend} (in a more general setting).

The required notation reads as follows:
 if $p$ and $m$ are positive integers, $0<p\leq m$,
$$
Q_{p,m}=\{[i_1, \dots, i_p]: 1\leq i_1<\dots <i_p\leq m, \; i_1, \ldots , i_p \in \mathbb{Z}\}.
$$
For $P(s)\in \FF[s]^{m \times n}$, $I\in Q_{p, m}$ and $J\in Q_{q, n}$, we denote by  $P(s)(I, J)$  the
$p\times q$ submatrix of $P(s)$ of entries that lie in the rows indexed by $I$ and the columns indexed by $J$.

Given $I=[i_1, \dots, i_p]\in Q_{p,m}$, let $I^*=[m-i_p+1, \dots, m-i_1+1]\in Q_{p,m}$, and if
$I'=[i'_1, \dots, i'_p]\in Q_{p,m}$, we write $I'\leq I$ whenever $i'_j\leq i_j$, $1\leq j \leq p$.

\medskip

The following lemma  is an immediate consequence of Lemma \ref{lemmamingID}.

\begin{lemma}\label{lemmaming}
Let $  E(s)=[e_{i,j}(s)]_{1\leq i,j\leq r}\in \FF[s]^{r \times r}$ with $\rank (  E(s))=r$. Let $k \in \{1, \dots, r\}$ and $Z\in  Q_{k,r}$. Then,
there exist $I, J\in Q_{k,r}$, $J\leq Z$,  $I\leq Z^*$, such that
$\det(  E(s)(I,J))\neq 0$.
\end{lemma}

{\bf Proof of the necessity of Theorem \ref{theomainlk2}.}
Let  $K(s)\in \FF[s]^{m\times r}$,  $L(s)^T\in \FF[s]^{n\times r}$ be minimal  bases
with column degrees $k_1\geq \dots\geq  k_r$ and $\ell_1\geq \dots\geq  \ell_r$, respectively.
Assume that there exists $A(s)\in\FF[s]^{m\times n}$ of  $\rank(A(s))=r$, $\deg(A(s))=d$, with
$\alpha_1(s)\mid \dots \mid \alpha_r(s)$ as invariant factors,   $f_1\leq  \dots\leq f_r$ as partial multiplicities of $\infty$, $\mR(A(s))=\mR(K(s))$, and $\mR(A(s)^T)=\mR(L(s)^T)$.
Let $g_1\geq \dots \geq g_r$ be the decreasing reordering of $\ell_1+k_{r}, \dots, \ell_r+k_1$.

 By Corollary \ref{corxy}, (\ref{eqf1}) holds.
It  only remains  to prove (\ref{eqprec}).
By Theorem \ref{theoreformulation}, there exists $  E(s)\in \FF[s]^{r\times r}$ with  invariant factors $\alpha_1(s)$, \ldots, $\alpha_r(s)$, such that the matrix $F(s)$ defined as in (\ref{eqmatF}) has degree $d$ and
$f_1,\dots, f_r$ as partial multiplicities of $\infty$.
Then, $F(s)$ satisfies  condition  (\ref{eqmaxminP2}) substituting $P(s)$ by $F(s)$, and where $M_k (F(s))$ is defined as in \eqref{def.MkP}.

Take $k\in \{1, \dots r-1\}$ and $Z\in Q_{k,r}$.
By Lemma \ref{lemmaming},
there exist $I, J\in Q_{k,r}$ such that $J\leq Z$,  $I\leq Z^*$ and
$\det(  E(s)(I,J))\neq 0$. Then,
$$
\det (F(s)(I,J))=s^{\sum_{j\in J}\ell_j+\sum_{i\in I }k_i}\det(  E(s)(I,J)),
$$
hence
$$
\deg(\det (F(s)(I,J)))= \sum_{j\in J}\ell_j+\sum_{i\in I}k_i+\deg(\det(  E(s)(I,J))).
$$
As
$M_k(F(s))\geq \deg(\det (F(s)(I,J)))$, $ \sum_{j\in J}\ell_j+\sum_{i\in I}k_i\geq \sum_{j\in Z}\ell_j+\sum_{i\in Z^*}k_i$, and $
\alpha_1(s)\cdots \alpha_k(s)\mid \det(  E(s)(I,J))\neq 0
$ (see \eqref{minork}),
we obtain
$$
M_k (F(s))-\sum_{i=1}^k\deg(\alpha_i)\geq \sum_{j\in Z}\ell_j+\sum_{i\in Z^*}k_i,
$$
therefore,
$$
M_k (F(s))-\sum_{i=1}^k\deg(\alpha_i)\geq
\max\{\sum_{j\in Z}\ell_j+\sum_{i\in Z^*}k_i\; : \; Z \in Q_{k,r}\}=\sum_{i=1}^kg_i.$$
From (\ref{eqmaxminP2})
we obtain that
$kd-\sum_{i=1}^kf_i-\sum_{i=1}^k\deg(\alpha_i)\geq \sum_{i=1}^kg_i$,
hence,
\begin{equation}\label{dggeq}
  \sum_{i=1}^k(d-g_i)\geq \sum_{i=1}^k\deg(\alpha_i)+\sum_{i=1}^kf_i, \quad 1\leq k \leq r-1.
\end{equation}
By Corollary \ref{corxy}, we get
\begin{equation}\label{dgeq}
  \sum_{i=1}^r(d-g_i)= \sum_{i=1}^r\deg(\alpha_i)+\sum_{i=1}^rf_i.
\end{equation}
Conditions (\ref{dggeq}) and (\ref{dgeq}) are equivalent to
(\ref{eqprec}).

\hfill $\Box$

\subsection{Proof of Theorem \ref{theomainlk2}: sufficiency of the conditions}
\label{subsecsuf}
The proof has two parts. First, in Section \ref{subsubf1}, we prove the sufficiency when the prescribed partial multiplicities at $\infty$ are all equal to zero, or, equivalently, when it is prescribed that infinity is not an eigenvalue of the polynomial matrix whose existence we wish to establish. Second, in Section \ref{subsubgen}, we prove the sufficiency in general. To achieve it, we use an argument based on M\"obius transformations (see, \cite[Section 6]{AmMaZa15} and \cite{MMMM15}). A M\"obius transformation is applied to the prescribed data to remove infinity as an eigenvalue. The transformed problem is solved by using the main result in Section \ref{subsubf1} ( i.e., Proposition \ref{propsufnoteid}). Finally, the polynomial matrix obtained in the previous step is transformed back with the inverse of the M\"obius transformation used in the first step. Strategies that use M\"obius transformations to avoid the eigenvalues at infinity have been used for solving other problems related to the existence of polynomial matrices with prescribed properties. See, for instance, \cite{taslaman-tiss-zab,tisseur-zaballa} and \cite{DeDoVa15}.

\subsubsection{Sufficiency when the prescribed partial multiplicities of $\infty$ are equal to $0$}\label{subsubf1}

The proof of the sufficiency in the case $f_1 = \cdots = f_r = 0$ follows from combining Theorem \ref{theoreformulation} with previous results on the existence of triangular regular polynomial matrices with prescribed invariant factors and diagonal entries. We recall first such previous results.

\begin{theorem} {\rm\cite[p. 208]{Sa80}} \label{lemasa2} Let
$\alpha_1(s) \mid \cdots \mid \alpha_r(s)$ and $\delta_1(s),\ldots,\delta_r(s)$
be monic polynomials with coefficients in an arbitrary field $\FF$. Then, there exists an $r\times r$ triangular
polynomial matrix with diagonal $(\delta_1(s),\ldots,\delta_r(s))$ and
$\alpha_1(s),\dots,\alpha_r(s)$ as invariant factors if and only if
\begin{gather} \label{alphadelta}
\alpha_1(s) \cdots \alpha_k(s) \mid \gcd \{\delta_{i_1}(s) \cdots  \delta_{i_k}(s) :
1\leq i_1<\dots < i_k\leq r\},\quad  1\leq k \leq r-1, \vspace{0.2cm}\\
\alpha_1(s)\cdots \alpha_r(s)=\delta_1(s)\cdots \delta_r(s). \label{alphadelta1}
\end{gather}
\end{theorem}

We will also need the following lemma, whose proof can be
found within that of \cite[Corollary 4.3]{Za97}. A result more general than Lemma \ref{lemaconher2}, together with a detailed proof, can be found in \cite[Lemma 4.6]{taslaman-tiss-zab}.
\begin{lemma}\label{lemaconher2}
  {\rm \cite[Lemma 4.3]{BaZa02}}
  Let $\FF$ be algebraically closed and let
$h_1,\ldots, h_r$ and $\alpha_1(s)\mid\cdots\mid\alpha_r(s)$ be non negative integers and
  monic polynomials with coefficients in $\FF$, respectively.
Let $g_1\geq \dots \geq g_r$ be the decreasing reordering of $h_1,\ldots, h_r$.
  If
\[
(g_1,\ldots, g_r)\prec (\deg(\alpha_r),\ldots, \deg(\alpha_1)),
\]
then, there exist monic polynomials $\delta_1(s),\ldots,\delta_r(s)$ such that $\deg(\delta_i)=h_i$,
$1\leq i\leq r$, and \eqref{alphadelta} and \eqref{alphadelta1} hold.
\end{lemma}

In the next proposition we prove the sufficiency part of Theorem \ref{theomainlk2} when we prescribe $f_1=\dots =f_r=0$.

\begin{proposition}\label{propsufnoteid}
Let $\FF$ be algebraically closed, $m$, $n$, $r\leq \min\{m,n \}$ be positive integers and $d$ a non negative integer.
Let
$\alpha_1(s)\mid \dots \mid \alpha_r(s)$ be monic polynomials  in $\FF [s]$ and let $K(s)\in \FF[s]^{m\times r}$,
$L(s)^T\in \FF[s]^{n\times r}$ be minimal bases with column degrees $k_1\geq \dots \geq  k_r$ and $\ell_1\geq \dots\geq  \ell_r$, respectively.
Let
$g_1\geq \dots \geq g_r$ be the decreasing  reordering of $k_r+\ell_1, \dots, k_1+\ell_r$.
If
\begin{equation*}\label{majnodei}
  (d-g_r,\dots, d-g_1)\prec (\deg(\alpha_r),\dots,\deg(\alpha_1)),
\end{equation*}
then, there exists a  polynomial matrix $A(s)\in\FF[s]^{m\times n}$ with $\rank(A(s))=r$, $\deg(A(s))=d$,  $\alpha_1(s), \dots, \alpha_r(s)$ as invariant factors,  $f_1=\dots= f_r=0$ as partial multiplicities of $\infty$,
$\mR(A(s))=\mR(K(s))$ and $\mR(A(s)^T)=\mR(L(s)^T)$.
\end{proposition}

{\bf Proof.}
By Lemma \ref{lemaconher2} and Theorem \ref{lemasa2}, there exist
monic polynomials $e_{1,1}(s),\ldots, \allowbreak e_{r,r}(s)$ such that
$$\deg(e_{i,i})=d-(k_{r-i+1}+\ell_i), \quad 1\leq i \leq r,$$
and an upper triangular matrix
$E'(s)\in \FF[s]^{r\times r}$
with diagonal ($e_{1,1}(s),\ldots,  e_{r,r}(s)$) and
$\alpha_1(s),\dots,\alpha_r(s)$ as invariant factors.
Next, inspired by \cite[Lemma 2.4]{tisseur-zaballa}, we prove that there exists a unimodular upper triangular polynomial matrix $U(s)\in \FF[s]^{r\times r}$ with diagonal entries all equal to one and such that
$$
E(s)=E'(s)U(s)=[e_{i,j}(s)]_{1\leq i,j\leq r}$$
satisfies
$$\deg(e_{i,j})<d-(k_{r-i+1}+\ell_i), \quad 1\leq i<j \leq r.
$$
The matrix $U(s)$ can be constructed by applying first to $E'(s)$ a sequence of elementary column replacement operations described with two nested ``for'' loops as follows: for $i = r-1, r-2, \ldots, 1$ and for $j=i+1, i+2, \ldots r$, replace the $j$th column by the $j$-th column minus the $i$-th column times the Euclidean quotient of the $(i,j)$-entry divided by the $(i,i)$-entry. Then, $U(s)$ is obtained by applying to $I_r$ the same sequence of elementary column replacement operations that has been applied to $E'(s)$.

Let
$$
F(s)=\diag(s^{k_r}, \dots, s^{k_1})  E(s)\diag(s^{\ell_1}, \dots, s^{\ell_r}).
$$
Then, $F(s)$ is upper triangular and if $F(s)=[f_{i,j}(s)]_{1\leq i,j\leq r}$,
$$
\begin{array}{ll}
  \deg(f_{i,i})=d, & 1\leq i \leq r,\\
\deg(f_{i,j})<d-(\ell_i+k_{r-i+1})+(\ell_j+k_{r-i+1})\leq d , &1\leq i<j \leq r.\\
  \end{array}
$$
As a consequence, $\deg(F(s))=d$ and, according to \eqref{def.MkP},
$$M_{k}(F(s))=kd, \quad 1\leq k \leq r.$$
Therefore, by \eqref{eqmaxminP2}, the partial multiplicities of $\infty$ of $F(s)$ are $f_1=\dots= f_r=0$. The proposition follows from  Theorem \ref{theoreformulation}  and Remark
\ref{remreformulationord}.
\hfill $\Box$

\subsubsection{Sufficiency in the general case}\label{subsubgen}

Our aim in this section is to prove the sufficiency of Theorem \ref{theomainlk2} in the general case, that is, when there may be non zero partial multiplicities of $\infty$.

{\bf Proof of the sufficiency of Theorem \ref{theomainlk2}.}
Assume that $\FF$ is an algebraically closed field and that conditions  (\ref{eqf1})
and (\ref{eqprec}) hold.
Let $a\in \efe$ be a  scalar that is not a root of $\alpha_r(s)$, i.e., $\gcd(\alpha_i(s),s-a)=1$ for $i=1,\dots,r$.
Let
$$\overline{K}(s)=K\left(\frac{1}{s}+a\right)\diag(s^{k_1},\dots,s^{k_r}),\  \overline{L}(s)^T=L\left(\frac{1}{s}+a\right)^T\diag(s^{\ell_1},\dots,s^{\ell_r}).
$$
By \cite[Theorem 7.4]{MMMM15}, $\overline{K}(s)$ and $\overline{L}(s)^T$ are minimal bases with column degrees $k_1,\dots,k_r$ and $\ell_1,\dots,\ell_r$, respectively. This can also be proved directly by combining Theorem \ref{mincharact} with some algebraic manipulations.

Write $\alpha_i(s)$ in the basis $\{1, s-a, (s-a)^2 , \ldots \}$, say,
$$
\alpha_i(s)=(s-a)^{\deg(\alpha_i)}+\dots + \alpha_{i,1} (s-a) + \alpha_{i,0},  \quad 1\leq i \leq r.
$$
Notice that $\alpha_{i,0}\neq 0$, $1\leq i\leq r$.
Then
\begin{equation}\label{eqalfatilde}
\tilde \alpha_i(s)=s^{\deg(\alpha_i)}\alpha_i\left(\frac{1}{s}+a\right)=1+\dots + \alpha_{i,1}s^{\deg(\alpha_i)-1} +\alpha_{i,0}s^{\deg(\alpha_i)},\quad 1\leq i \leq r,
\end{equation}
is a polynomial of the same degree as $\alpha_i(s)$ and such that $\gcd(\tilde \alpha_i(s),s)=1$. Define
\begin{equation}\label{eqbeta}
\beta_i(s)=\frac{1}{\alpha_{i,0}} \tilde \alpha_i(s)s^{f_i},\quad 1\leq i \leq r.
\end{equation}
These are monic polynomials of degree $\deg(\alpha_i)+f_i$ and $\beta_1(s)\mid\dots\mid\beta_r(s)$.

Condition (\ref{eqprec}) is equivalent to
$$ (d-g_r,\dots, d-g_1)\prec (\deg(\beta_r),\dots,\deg(\beta_1)).$$
By Proposition \ref{propsufnoteid} there exists $B(s)\in\FF[s]^{m\times n}$, with $\rank(B(s))=r$, $\deg(B(s))=d$,
$\beta_1(s),\dots, \beta_r(s)$ as invariant factors,  $e_1=\dots=e_r=0$ as partial multiplicities of $\infty$, $\mR(B(s))=\mR(\overline{K}(s))$ and $\mR(B(s)^T)=\mR(\overline{L}(s)^T)$.

Assume that
$
B(s)=B_ds^d+ B_{d-1} s^{d-1} + \dots+B_0.
$
Notice that $B_0\neq 0$, because $\gcd(\tilde \alpha_1(s),s)=1$, condition (\ref{eqf1}) holds, and, therefore, $\gcd(\beta_1(s),s)=1$.
Set
\begin{equation} \label{eq.defAs}
A(s)=(s-a)^d B\left(\frac{1}{s-a}\right)=B_d+ B_{d-1} (s-a) + \dots+B_0(s-a)^d.
\end{equation}
It is clear that $A(s)\in\efe[s]^{m\times n}$, $\rank(A(s))=r$, and $\deg(A(s))=d$. Furthermore, $$\overline{K}\left(\frac{1}{s-a}\right)\diag\left((s-a)^{k_1},\dots,(s-a)^{k_r}\right)=K(s),
$$
$$
\overline{L}\left(\frac{1}{s-a}\right)^T\diag\left((s-a)^{\ell_1},\dots,(s-a)^{\ell_r}\right)=L(s)^T.
$$
Next, observe that $\mR (A(s)) = \mR (B(\frac{1}{s-a}))$ because $(s-a)^d$ is a scalar in $\FF(s)$. Moreover, $B(s) = \overline{K} (s) C(s)$, where $C(s) \in \FF[s]^{r \times n}$ is the matrix whose $i$-th column contains the coordinates of the $i$-th column of $B(s)$ in the minimal basis $\overline{K} (s)$ for $1 \leq i \leq n$ ($C(s)$ is polynomial by \cite[p. 495]{Fo75}). Thus, $B(\frac{1}{s-a}) = \overline{K} (\frac{1}{s-a}) C(\frac{1}{s-a})$, which implies that $\mR (A(s)) = \mR (B(\frac{1}{s-a})) \subseteq \mR (\overline{K} (\frac{1}{s-a})) = \mR (K(s))$. The previous inclusion combined with $\dim \mR (A(s)) = \dim \mR (K(s)) = r$ implies that $\mR(A(s))=\mR(K(s))$. A similar argument proves that
$\mR(A(s)^T)=\mR(L(s)^T)$.

Let $b_1,\dots,b_h$ be the non zero finite eigenvalues of $B(s)$.
Write
\begin{equation}\label{eqbetafi}
\beta_i(s)=\prod_{j=1}^{h}(s-b_j)^{t_{ij}}s^{f_i},\quad 1\leq i\leq r.
\end{equation}

By taking $\alpha=0, \beta=1, \gamma=1,$ and $\delta=-a$ in \cite[Proposition 6.16]{AmMaZa15} (see also \cite[Theorem 5.3]{MMMM15}) we deduce:
\begin{itemize}
\item [(i)] For any $c\in\efe\setminus\{0\}$,
if $(s-c)^{t_{1}},\dots,(s-c)^{t_{r}}$ are the local invariant rational functions at $s-c$ of $B(s)$ then $\left(s-a-\frac{1}{c}\right)^{t_1},\dots,\left(s-a-\frac{1}{c}\right)^{t_r}$ are the local invariant rational functions at $s-a-\frac{1}{c}$ of $A(s)$. More precisely, $\left(s-a-\frac{1}{b_j}\right)^{t_{1j}},\dots, \allowbreak \left(s-a-\frac{1}{b_j}\right)^{t_{rj}}$ are the local invariant rational functions at $s-a-\frac{1}{b_j}$ of $A(s)$, for $j=1,\dots,h$, and for any other scalar $c'\neq a$, i.e., $c'\in\efe\setminus\{a,a + \frac{1}{b_1},\dots, a+ \frac{1}{b_h}\}$ the local invariant rational functions at $s-c'$ of $A(s)$ are $1,\dots,1$.

\item [(ii)] As $f_1,\dots,f_r$ are the partial multiplicities of $0$ in $B(s)$, the partial multiplicities of $\infty$ in $A(s)$ are $0, f_2-f_1,\dots,f_r-f_1$. As condition (\ref{eqf1}) holds, the partial multiplicities of $\infty$ in $A(s)$ are $f_1,\dots, f_r$.

\item [(iii)] Since $e_1=\dots=e_r=0$ are the partial multiplicities of $\infty$ in $B(s)$, the local invariant rational functions at $s-a$ of $A(s)$ are $1,\dots, 1$, which means that the partial multiplicities of $a$ in $A(s)$ are $0,\dots,0$.
\end{itemize}

By (i) and (iii) the invariant factors of $A(s)$ are $\prod_{j=1}^{h}\left(s-a-\frac{1}{b_j}\right)^{t_{ij}}, 1\leq i \leq r$. It only remains to prove that they are equal to $\alpha_1 (s), \ldots, \alpha_r (s)$.

From (\ref{eqbeta}) and (\ref{eqbetafi}) we obtain
$\frac{1}{\alpha_{i,0}} \tilde \alpha_i(s)=
\prod_{j=1}^{h}(s-b_j)^{t_{ij}}$ and, therefore, $\deg(\tilde{\alpha}_i)=\sum_{j=1}^h t_{ij}$ and $1=\tilde\alpha_i(0)=\alpha_{i,0}\prod_{j=1}^{h}(-b_j)^{t_{ij}}, 1\leq i \leq r$.
Now, from (\ref{eqalfatilde}),
$$\tilde\alpha_i\left(\frac{1}{s-a}\right)=
\left(\frac{1}{s-a}\right)^{\deg(\alpha_i)}\alpha_i(s), \quad 1\leq i\leq r.
$$
Thus,
$$
\begin{array}{ll}
\alpha_i(s)&=(s-a)^{\deg(\alpha_i)}\tilde\alpha_i\left(\frac{1}{s-a}\right)=(s-a)^{\deg(\alpha_i)}\alpha_{i,0}\prod_{j=1}^{h}\left(\frac{1}{s-a}-b_j\right)^{t_{ij}}\\
&=\alpha_{i,0}(s-a)^{\deg(\alpha_i)}\prod_{j=1}^{h}\left(\frac{-b_j}{s-a}\right)^{t_{ij}}\left(s-a-\frac{1}{b_j}\right)^{t_{ij}},\quad 1\leq i \leq r.
\end{array}
$$
Since $\deg(\alpha_i)=\deg(\tilde{\alpha}_i)=\sum_{j=1}^h t_{ij}$ and $1=\alpha_{i,0}\prod_{j=1}^{h}(-b_j)^{t_{ij}}$,
$$
\alpha_i(s)=\prod_{j=1}^{h}\left(s-a-\frac{1}{b_j}\right)^{t_{ij}},\quad 1\leq i \leq r.
$$

Summing up, the polynomial matrix $A(s)$ defined in \eqref{eq.defAs} satisfies $A(s)\in\efe[s]^{m\times n}$, $\rank(A(s))=r$, $\deg(A(s))=d$,  $\alpha_1(s),\dots,\alpha_r(s)$ are its invariant factors, $f_1,\dots,f_r$ are its partial multiplicities of $\infty$ (see (ii)), $\mR(A(s))=\mR(K(s))$ and $\mR(A(s)^T)=\mR(L(s)^T)$.
\hfill $\Box$

\section{Results for rational matrices} \label{sec.rational}

In this section, we solve problems for rational matrices similar to those stated for polynomial matrices in Section \ref{problems}. In Theorem \ref{theomainlk2_rat}, we solve the rational counterpart of Problem \ref{probkl2} by combining the polynomial result in Theorem \ref{theomainlk2} with the simple  fact that if $q(s)$ is the least common denominator of the entries of a rational matrix $R(s)$, then $q(s) R(s)$ is a polynomial matrix. It allows us to ``transform'' the rational problem of prescribed data into a polynomial problem with ``related'' prescribed data. We then solve the resulting polynomial problem using Theorem \ref{theomainlk2}, and finally transform its solution into the solution of the original rational problem. The properties of the involved ``transformations'' rely on the auxiliary Lemma \ref{lem_polrat}. A similar strategy was used in \cite[Theorem 4.1]{anguasetal2019} to prove the rational counterpart of Theorem \ref{theoexistenceDeDoVa152}. The solutions of the rational counterparts of Problems \ref{probkl} and \ref{probklcv} follow easily from Theorem \ref{theomainlk2_rat} and are given in Theorems \ref{cormainlk_rat} and \ref{cormainlkcv_rat}, respectively.

In Lemma \ref{lem_polrat}, we relate the finite and infinite invariant rational functions, the right and left null spaces, and the column and row spaces of a rational matrix with those of a polynomial matrix related to it. In the rest of this section, it is convenient to bear in mind that the denominator $\psi_1 (s)$ of the first invariant rational function of a rational matrix $R(s)$ is precisely the (monic) least common denominator of the entries of $R(s)$. See \cite[p. 444]{Kail80}.

\begin{lemma}\label{lem_polrat}
Let $R(s) \in \FF(s)^{m\times n}$ be a rational matrix and let $p(s) \in \FF(s)$ be a monic polynomial multiple of the least common denominator of the entries of $R(s)$. Then, $p(s)R(s) \in \FF[s]^{m\times n}$ is a polynomial matrix with the same rank as $R(s)$ and:
\begin{itemize}
	\item[\rm(i)] $\frac{\epsilon_1(s)}{\psi_1(s)},\dots,\frac{\epsilon_r(s)}{\psi_r(s)}$ are the invariant rational functions of $R(s)$ if and only if $\frac{p(s)\epsilon_1(s)}{\psi_1(s)},\dots, \allowbreak \frac{p(s)\epsilon_r(s)}{\psi_r(s)}$ are the invariant factors of $p(s)R(s)$.
	\item[\rm(ii)] The integers $q_1,\dots,q_r$ are the invariant orders at infinity of $R(s)$ if and only if $\deg(p(s)R(s))=\deg(p(s))-q_1$ and  $0=q_1-q_1, q_2-q_1,\dots,q_r-q_1$ are the partial multiplicities of $\infty$ in $p(s)R(s)$.
	\item[\rm(iii)] $\mathcal{N}_r(R(s))=\mathcal{N}_r(p(s)R(s))$, $\mathcal{N}_{\ell}(R(s))=\mathcal{N}_{\ell}(p(s)R(s))$, and therefore the right (left) minimal indices of $p(s)R(s)$ and of $R(s)$ are equal.
	\item[\rm(iv)] $\mR(p(s)R(s))=\mR(R(s))$, $\mR((p(s)R(s))^T)=\mR(R(s)^T)$, and therefore the col-span (row-span) minimal indices of $p(s)R(s)$ and of $R(s)$ are equal.
\end{itemize}
\end{lemma}
{\bf Proof.} The proof is straightforward. So, we just provide some hints and prove with some detail only item (ii). It is obvious that $p(s) R(s)$ is a polynomial matrix. The equality $\rank(R(s)) = \rank (p(s) R(s))$ and items (iii) and (iv) follow from the facts that $p(s)$ is a scalar in $\FF(s)$ and that multiplication by non zero scalars in the field does not change the null spaces and the column and row spaces of a matrix. Item (i) follows from the Smith--McMillan form \eqref{eq.smithmcmillan}. To prove (ii), note that $p(s) R(s) =  \left(\frac{p(s)}{s^{\deg(p)}} I_m \right) (s^{\deg(p)} R(s))$ and that the matrix  $\left(\frac{p(s)}{s^{\deg(p)}} I_m \right)$ is biproper. Thus, the Smith-McMillan form at infinity of $p(s) R(s)$ is equal to the Smith-McMillan form at infinity of  $s^{\deg(p)} R(s)$, which is the Smith-McMillan form at infinity of $R(s)$ multiplied by $s^{\deg(p)}$. Therefore, taking into account \eqref{eq.infsmithmcmillan}, $q_1,\dots,q_r$ are the invariant orders at infinity of $R(s)$ if and only if $q_1 - \deg(p) ,\dots,q_r - \deg(p)$ are the invariant orders at infinity of $p(s) R(s)$. Then $\deg(p(s)R(s))=\deg(p(s))-q_1$ and, by \eqref{eqqf},  $0=q_1-q_1, q_2-q_1,\dots,q_r-q_1$ are the partial multiplicities of $\infty$ in $p(s)R(s)$. The converse also follows from \eqref{eqqf}.
\hfill $\Box$

\begin{theorem}\label{theomainlk2_rat}
	Let $m$, $n$, $r\leq \min\{m,n \}$ be  positive integers.
	Let
	$\epsilon_1(s)\mid \dots \mid \epsilon_r(s)$ and $\psi_r(s)\mid \dots \mid \psi_1(s)$ be monic polynomials  in $\FF[s]$ such that $\frac{\epsilon_i(s)}{\psi_i(s)}$ are irreducible rational functions for $i=1,\dots, r$. Let $q_1\leq\dots\leq q_r$ be integers  and  let $K(s)\in \FF[s]^{m\times r}$ and  $L(s)^T\in \FF[s]^{n\times r}$  be minimal  bases
	with column degrees $k_1\geq \dots\geq  k_r$ and $\ell_1\geq \dots\geq  \ell_r$, respectively.
	Let
	$g_1\geq \dots \geq g_r$ be the decreasing reordering of $k_r+\ell_1, \dots, k_1+\ell_r$.
	
	If $R(s)\in\FF(s)^{m\times n}$ is a rational matrix with $\rank(R(s))=r$,
	$\frac{\epsilon_1(s)}{\psi_1(s)},\dots,\frac{\epsilon_r(s)}{\psi_r(s)}$ as invariant rational functions, $q_1, \dots, q_r$ as invariant orders at $\infty$, $\mR(R(s))=\mR(K(s))$ and $\mR(R(s)^T)=\mR(L(s)^T)$,  then
 	\begin{equation}\label{eqprec_rat1}
			(-g_r, \dots, -g_1)\prec (\deg(\epsilon_r)-\deg(\psi_r)+q_r, \dots, \deg(\epsilon_1)-\deg(\psi_1)+q_1).
		\end{equation}
	
	Conversely, if $\FF$ is  algebraically closed,  condition \eqref{eqprec_rat1} is  sufficient for the existence of a rational matrix $R(s)\in\FF(s)^{m\times n}$ with $\rank(R(s))=r$,
	$\frac{\epsilon_1(s)}{\psi_1(s)},\dots,\frac{\epsilon_r(s)}{\psi_r(s)}$ as invariant rational functions, $q_1, \dots, q_r$ as invariant orders at $\infty$, $\mR(R(s))=\mR(K(s))$ and $\mR(R(s)^T)=\mR(L(s)^T)$.
\end{theorem}

\begin{rem}\label{remprec_rat} \textit{ (Comparing Theorems \ref{theomainlk2_rat} and \ref{theomainlk2}).}
{\rm

\begin{enumerate}
\item
		Clearly, when the rational matrix $R(s)$ is polynomial, the majorization (\ref{eqprec_rat1}) reduces to condition (\ref{eqprec}).
  Namely, $\psi_i(s)=1$  for $i=1,\dots,r$, $\epsilon_1(s),\dots,\epsilon_r(s)$ are the invariant factors of $R(s)$, $\deg(R(s))=-q_1$ and the partial multiplicities of $\infty$ in $R(s)$ are $q_1-q_1,\dots,q_r-q_1$ by \eqref{eqqf}.

 \item  Observe also that  (\ref{eqprec_rat1}) is equivalent to
\begin{equation}\label{eqcondleq_rat}
		\sum_{i=1}^k g_i+\sum_{i=1}^k\deg(\epsilon_i)-\sum_{i=1}^k\deg(\psi_i)+\sum_{i=1}^k q_i\leq 0, \quad 1\leq k\leq r-1,
	\end{equation}
	\begin{equation*}\label{eqsumklfa_rat}
		\sum_{i=1}^r k_i+\sum_{i=1}^r \ell_i+\sum_{i=1}^r\deg(\epsilon_i)-\sum_{i=1}^r\deg(\psi_i)+\sum_{i=1}^r q_i=0.
	\end{equation*}
    Taking into account \eqref{eqsums}, the above equality is the index sum theorem constraint for rational matrices (see \cite[Remark 3.4]{anguasetal2019} and \cite[Theorem 2.4]{AmBaMaRo25}). Altogether, this remark can be seen as the ``rational'' counterpart of Remark \ref{remprec}-3.
    \item It is interesting to point out that the only condition involved in Theorem \ref{theomainlk2_rat} is \eqref{eqprec_rat1} in contrast with the two conditions \eqref{eqf1} and \eqref{eqprec} involved in Theorem \ref{theomainlk2}. The reason of this difference is that there exist rational matrices with arbitary values of the smallest invariant order at infinity, while $f_1 = 0$ always holds for polynomial matrices. \hfill $\Box$
 \end{enumerate}
 }
\end{rem}

{\bf Proof of Theorem \ref{theomainlk2_rat}}.
Assume that there exists a rational matrix $R(s)\in\FF(s)^{m\times n}$ with $\rank(R(s))=r$,
$\frac{\epsilon_1(s)}{\psi_1(s)},\dots,\frac{\epsilon_r(s)}{\psi_r(s)}$ as invariant rational functions, $q_1, \dots, q_r$ as invariant orders at $\infty$, $\mR(R(s))=\mR(K(s))$ and $\mR(R(s)^T)=\mR(L(s)^T)$.  By Lemma \ref{lem_polrat}, $\psi_1(s)R(s)$ is a polynomial matrix with  $\rank(\psi_1(s)R(s))=r$, $\deg(\psi_1(s)R(s))=\deg(\psi_1(s))-q_1$, $\frac{\psi_1(s)\epsilon_1(s)}{\psi_1(s)},\dots,$ $\frac{\psi_1(s)\epsilon_r(s)}{\psi_r(s)}$ as invariant factors, $q_1-q_1,q_2-q_1,\dots,q_r-q_1$ as partial multiplicities of $\infty$, $\mR(\psi_1(s)R(s))=\mR(R(s))=\mR(K(s))$ and $\mR((\psi_1(s)R(s))^T)=\mR(R(s)^T)=\mR(L(s)^T)$. By Theorem \ref{theomainlk2} applied to $\psi_1 (s) R(s)$, we have
\begin{equation} \label{eqprec_rat2}\begin{array}{l}(\deg(\psi_1)-q_1-g_r,\dots, \deg(\psi_1)-q_1-g_1)  \\ \phantom{aaaaaa} \prec
(\deg(\psi_1)+\deg(\epsilon_r)-\deg(\psi_r)+q_r-q_1,\dots, \\ \phantom{aaaaaa \prec
(} \deg(\psi_1)+\deg(\epsilon_1)-\deg(\psi_1)+q_1-q_1),
\end{array}
\end{equation}
which is equivalent to condition (\ref{eqprec_rat1}).

Conversely, assume that $\FF$ is  algebraically closed and (\ref{eqprec_rat1})  holds. From (\ref{eqcondleq_rat}) for $k=1$, we obtain $g_1+\deg(\epsilon_1)\leq \deg(\psi_1)-q_1$, therefore, $\deg(\psi_1)-q_1\geq 0$. Moreover,
notice that $\alpha_i(s)=\frac{\psi_1(s)\epsilon_i(s)}{\psi_i(s)}$, $i=1\dots,r$, are monic polynomials such that $\alpha_1(s)\mid\dots\mid\alpha_r(s)$ and $\deg (\alpha_i) = \deg (\psi_1) + \deg (\epsilon_i) - \deg(\psi_i)$. Note also that $(q_r-q_1,\dots,q_1-q_1)$ is a partition of non negative integers. As already seen in the necessity part, condition (\ref{eqprec_rat1}) is equivalent to (\ref{eqprec_rat2}). By Theorem \ref{theomainlk2}, there exists $A(s)\in\FF[s]^{m\times n}$
with $\rank(A(s))=r$, $\deg(A(s))=\deg(\psi_1)-q_1$,
$\alpha_1(s),\dots, \alpha_r(s)$ as invariant factors,  $q_1-q_1, \dots, q_r-q_1$ as partial multiplicities of $\infty$, $\mR(A(s))=\mR(K(s))$ and $\mR(A(s)^T)=\mR(L(s)^T)$. Take $R(s)=\frac{A(s)}{\psi_1(s)}$. Then, $\rank(R(s))=\rank(A(s))=r$, and by Lemma \ref{lem_polrat}, $\frac{\epsilon_1(s)}{\psi_1(s)},\dots,\frac{\epsilon_r(s)}{\psi_r(s)}$ are the invariant rational functions of $R(s)$, $q_1, \dots, q_r$  its invariant orders at $\infty$, $\mR(R(s))=\mR(K(s))$ and $\mR(R(s)^T)=\mR(L(s)^T)$.
\hfill $\Box$

\begin{theorem}\label{cormainlk_rat}
Let $m$, $n$, $r\leq \min\{m,n \}$ be  positive integers. Let $\epsilon_1(s)\mid \dots \mid \epsilon_r(s)$ and $\psi_r(s)\mid \dots \mid \psi_1(s)$ be monic polynomials in $\FF [s]$ such that $\frac{\epsilon_i(s)}{\psi_i(s)}$ are irreducible rational functions for $i=1,\dots, r$. Let $q_1\leq\dots\leq q_r$ be integers and  $(k_1, \ldots, k_{r})$, $(\ell_1, \ldots, \ell_{r})$ partitions of non negative integers. Let
$g_1\geq \dots \geq g_r$ be  the decreasing reordering of $k_r+\ell_1, \dots, k_1+\ell_r$.

If $R(s)\in\FF(s)^{m\times n}$ is a rational matrix with $\rank(R(s))=r$,
$\frac{\epsilon_1(s)}{\psi_1(s)},\dots,\frac{\epsilon_r(s)}{\psi_r(s)}$ as invariant rational functions, $q_1, \dots, q_r$ as invariant orders at $\infty$,
$k_1,\dots, k_r$ as col-span minimal indices, and  $\ell_1, \dots,\ell_r$ as row-span minimal indices, then \eqref{eqx>0}, \eqref{eqy>0} and \eqref{eqprec_rat1}  hold.

Conversely, if $\FF$ is  algebraically closed, conditions \eqref{eqx>0}, \eqref{eqy>0} and \eqref{eqprec_rat1} are sufficient for the existence of a rational matrix $R(s)\in\FF(s)^{m\times n}$ with $\rank(R(s))=r$,
$\frac{\epsilon_1(s)}{\psi_1(s)},\dots,\frac{\epsilon_r(s)}{\psi_r(s)}$ as invariant rational functions, $q_1, \dots, q_r$ as invariant orders at $\infty$, $k_1,\dots, k_r$ as col-span minimal indices and $\ell_1, \dots,\ell_r$ as row-span minimal indices.
\end{theorem}

{\bf Proof.} The  proof follows from Proposition \ref{exminbasis} and Theorem \ref{theomainlk2_rat} analogously as the proof of Theorem \ref{cormainlk}, i.e, the polynomial case, follows from Proposition \ref{exminbasis} and Theorem \ref{theomainlk2}.
\hfill $\Box$

\begin{theorem}\label{cormainlkcv_rat}
Let $m$, $n$, $r\leq \min\{m,n \}$ be  positive integers. Let
$\epsilon_1(s)\mid \dots \mid \epsilon_r(s)$ and $\psi_r(s)\mid \dots \mid \psi_1(s)$ be monic polynomials  in $\FF [s]$ such that $\frac{\epsilon_i(s)}{\psi_i(s)}$ are irreducible rational functions for $i=1,\dots, r$. Let $q_1\leq\dots\leq q_r$ be integers. Let
$(k_1, \ldots, k_{r})$, $(\ell_1, \ldots, \ell_{r})$, $(d_1, \dots, d_{n-r})$,  $(v_1, \dots, v_{m-r})$  be partitions of non negative integers. Let
$g_1\geq \dots \geq g_r$ be the decreasing  reordering of $k_r+\ell_1, \dots, k_1+\ell_r$.

If $R(s)\in\FF(s)^{m\times n}$ is a rational matrix with $\rank(R(s))=r$,
$\frac{\epsilon_1(s)}{\psi_1(s)},\dots,\frac{\epsilon_r(s)}{\psi_r(s)}$ as invariant rational functions, $q_1, \dots, q_r$ as invariant orders at $\infty$, $k_1,\dots, k_r$ as col-span minimal indices, $\ell_1, \dots,\ell_r$ as row-span minimal indices, and $d_1, \dots, d_{n-r}$ and $v_1, \dots, v_{m-r}$ as right and left minimal indices, respectively, then \eqref{eqsums}  and \eqref{eqprec_rat1} hold.

Conversely, if $\FF$ is  algebraically closed, conditions \eqref{eqsums} and \eqref{eqprec_rat1} are sufficient for the existence of a rational matrix $R(s)\in\FF(s)^{m\times n}$ with $\rank(R(s))=r$,
$\frac{\epsilon_1(s)}{\psi_1(s)},\dots,\frac{\epsilon_r(s)}{\psi_r(s)}$ as invariant rational functions, $q_1, \dots, q_r$ as invariant orders at $\infty$, $k_1,\dots, k_r$ as col-span minimal indices,  $\ell_1, \dots,\ell_r$ as row-span minimal indices, and
$d_1, \dots, d_{n-r}$ and  $v_1, \dots, v_{m-r}$ as right and left minimal indices, respectively.
\end{theorem}

{\bf Proof.}
The  proof is analogous to that of the polynomial case in Theorem \ref{cormainlkcv}. Simply, one has to replace the use of Theorems \ref{theomainlk2} and \ref{cormainlk} in that proof by Theorems \ref{theomainlk2_rat} and \ref{cormainlk_rat}, respectively.
\hfill $\Box$

\section{Conclusions and open problems} \label{sec.conclusions}
In this paper, we have provided necessary and sufficient conditions for the existence of polynomial and rational matrices when the minimal indices of the column and row spaces, together with other structural data, are prescribed. As far as we know, these are the first results available in the literature that provide complete information about how the minimal indices of the column and row spaces are related to other relevant magnitudes of rational and polynomial matrices, as the invariant rational functions, the invariant rational functions at infinity, and the minimal indices of the left and right null spaces.

It is well known that polynomial and rational matrices arising in applications often possess particular structures as, for instance, symmetric, skew-symmetric, Hermitian, palindromic or alternating structures \cite{MMMM-good}. Therefore, it would be of interest to study how the necessary and sufficient conditions obtained in this paper should be modified when the polynomial or rational matrix whose existence is to be established is required to have some particular structure. This is, in general, a very challenging problem, since extensions of the result in \cite{AmBaMaRo23,DeDoVa15} about the existence of polynomial matrices with prescribed complete eigenstructure (i.e., without considering the minimal indices of the column and row spaces) to structured scenarios have been obtained only in some cases and under some strong constraints as imposing regularity \cite{batzke-mehl-2014}, or simple eigenvalues \cite[Theorem 4.1]{DeTDmyDop2020},  or the degree to be at most two  \cite{quadratic-palin,perovic-mackey}. However, these constraints are significant in applications and may make it possible to obtain structured analogues of the results in this paper.

\begin{appendices}
\section{A determinantal lemma}\label{secappend}

Let ${\mathcal R}$ be an integral domain,
${\mathcal R}^{m\times n}$ the set of $m \times n$  matrices with entries in ${\mathcal R}$ and
${\mathcal F}$ the field of fractions of ${\mathcal R}$.
As in Section \ref{subsecnec}, for $P\in {\mathcal R}^{m \times n}$, $I\in Q_{p, m}$ and $J\in Q_{q, n}$, we denote by  $P(I, J)$  the
$p\times q$ submatrix of  $P$ of entries that lie in the rows indexed by $I$ and the
columns indexed by $J$, and $P(I, :)\in{\mathcal R}^{p\times n}$ and $P(:, J)\in {\mathcal R}^{m\times q}$ are the submatrices of $P$ formed by the rows indexed by $I$ and the columns indexed by $J$, respectively. Moreover, we also use the following notation: given two positive integers $i$ and $j$, such that $i < j$, $i:j = [i, i+1, \ldots, j]$.

\begin{lemma}\label{lemmamingID}
Let $E= [e_{i,j}]_{1\leq i,j\leq r}\in{\mathcal R}^{r \times r}$ with $\rank (E)=r$. Let $k \in \{1, \dots, r\}$ and $Z\in  Q_{k,r}$. Then,
there exist $I, J\in Q_{k,r}$, $J\leq Z$,  $I\leq Z^*$, such that
$\det(E(I,J))\neq 0$.
\end{lemma}

Lemma \ref{lemmamingID}  is proven by induction on $r$. First, we illustrate the proof with an example.

\begin{example}\label{ejproof}
{\rm
    Let $E= [e_{i,j}]_{1\leq i,j\leq 5}\in{\mathcal R}^{5 \times 5}$ with $\rank (E)=5$.

    Assume that if $\widehat E\in{\mathcal R}^{4 \times 4}$ has $\rank (\widehat E)=4$, then for $k\in \{1, \dots, 4\}$ and  $\hat Z\in Q_{k,4}$ there exist
 $\hat I, \hat J\in Q_{k,4}$, $\hat J\leq \hat Z$,  $\hat I\leq \hat Z^*$, such that
$\det(  \widehat E(\hat I,\hat J))\neq 0$.

    If
    $X=\begin{bmatrix}c_1&\cdots &c_5\end{bmatrix}\in{\mathcal R}^{4\times 5}$
    is
    the matrix formed by the $4$ first rows of $E$, then $\rank(X)=4$. Thus, at least one of the columns of $X$ depends linearly on the others.
Suppose that the columns $c_1, c_2$ are linearly independent and $c_3$ depends linearly on $c_1, c_2$.
Then, there exists $g\in {\mathcal F}^{2\times 1}$ such that $c_3=\begin{bmatrix}c_1&c_2\end{bmatrix}g$.
As $\rank (\begin{bmatrix}c_1&c_2&c_3\\e_{5,1}&e_{5,2}&e_{5,3}\end{bmatrix})=\rank (E(:, 1:3))=3$, we conclude that
$e_{5,3}\neq \begin{bmatrix}e_{5,1}&e_{5,2}\end{bmatrix}g$.

Take  $\widehat E= [\hat e_{i,j}]_{1\leq i,j\leq 4} =\begin{bmatrix}c_1&c_2&c_4 &c_5 \end{bmatrix}\in{\mathcal R}^{4\times 4}$;
i.e., $\widehat E$ is the submatrix obtained  from $E$  by deleting the last row and the third column. We have $\rank(\widehat E)=4$,
and
$$
\widehat E=\begin{bmatrix}
\hat e_{1,1}&\hat e_{1,2}&\hat e_{1,3}&\hat e_{1,4}\\\hat e_{2,1}&\hat e_{2,2}&\hat e_{2,3}&\hat e_{2,4}\\\hat e_{3,1}&\hat e_{3,2}&\hat e_{3,3}&\hat e_{3,4}\\\hat e_{4,1}&\hat e_{4,2}&\hat e_{4,3}&\hat e_{4,4}
\end{bmatrix}
=
\begin{bmatrix}
e_{1,1}&e_{1,2}&e_{1,4}&e_{1,5}\\e_{2,1}&e_{2,2}&e_{2,4}&e_{2,5}\\e_{3,1}&e_{3,2}&e_{3,4}&e_{3,5}\\e_{4,1}&e_{4,2}&e_{4,4}&e_{4,5}
\end{bmatrix}.
$$

Let $k \in \{1, \dots, 5\}$ and $Z=[z_1, \dots, z_k]\in  Q_{k,5}$.

If $k=5$, then $Q_{5,5}=\{[1, \dots, 5]\}$ and $Z=Z^*=[1, \dots, 5]$. Taking $I=J=Z$ we have $E(I,J)=E$ and $\det(E)\neq 0$.

Assume that $k \in \{1, \dots, 4\}$. We show how to proceed in two examples for $k=3$.
\begin{itemize}
\item Let $Z=[z_1, z_2, z_3]=[1, 3, 4]\in Q_{3, 5}$ ($Z^*=[2,3,5]$). Here we take
$\hat Z=[\hat z_1, \hat z_2, \hat z_3]=[1,2,3]\in  Q_{3, 4}$ ($\hat Z^*=[2, 3, 4]$). As we will see in the proof of Lemma \ref{lemmamingID}, the choice corresponds to $\hat z_1=z_1$, $\hat z_2=z_2-1$,  $\hat z_3=z_3-1$. By the induction hypothesis,
there exist
 $\hat I=[\hat i_1,\hat i_2 ,\hat i_3], \hat J=[\hat j_1,\hat j_2 ,\hat j_3]\in Q_{3,4}$ such that
$$\hat i_1\leq 2,\quad  \hat i_2\leq 3,\quad  \hat i_3 \leq 4, \quad \hat j_1=1, \quad \hat j_2= 2,\quad  \hat j_3= 3,$$ and $\det (\widehat E(\hat I, \hat J))\neq 0$.
If we take $I=\hat I$ and  $J=[1, 2, 4]$, then $I, J\in Q_{3,5}$,  $I\leq Z^*$, $J\leq Z$, and
$$
E(I,J)=
\begin{bmatrix}e_{i_1,1}&e_{i_1,2}&e_{i_1,4}\\e_{i_2,1}&e_{i_2,2}&e_{i_2,4}\\e_{i_3,1}&e_{i_3,2}&e_{i_3,4}\\\end{bmatrix}=
\begin{bmatrix}
\hat e_{\hat i_1,1}&\hat e_{\hat i_1,2}&\hat e_{\hat i_1,3}\\
\hat e_{\hat i_2,1}&\hat e_{\hat i_2,2}&\hat e_{\hat i_2,3}\\\hat e_{\hat i_3,1}&\hat e_{\hat i_3,2}&\hat e_{\hat i_3,3}
\end{bmatrix},
$$
i.e., $E(I,J)=\widehat E(\hat I, \hat J)$; hence $\det(E(I,J))\neq 0$.

\item Let $Z=[z_1,z_2,z_3]=[1, 2,3]\in Q_{3, 5}$ ($Z^*=[3,4,5]$).  Now we take
  $\hat Z=[\hat z_1, \hat z_2]=[1,2]\in  Q_{2, 4}$ ($\hat Z^*=[3,4]$). Again, as we will see in the proof of Lemma \ref{lemmamingID}, $\hat Z$ belongs  to  $Q_{k-1, 4}$ and is defined as $\hat z_1=z_1$,
  $\hat z_2=z_2$.
  By the induction hypothesis,
there exist
 $\hat I=[\hat i_1,\hat i_2], \hat J=[\hat j_1,\hat j_2]\in Q_{2,4}$ such that
$$\hat i_1\leq 3,\quad  \hat i_2\leq 4,\quad   \hat j_1=1, \quad \hat j_2= 2,$$ and $\det (\widehat E(\hat I, \hat J))\neq 0$.
Define $I=[\hat i_1, \hat i_2, 5]$ and $J=[1, 2, 3]$. Then $I, J\in Q_{3,5}$,  $I\leq Z^*$, $J\leq Z$, and
$$
E(I,J)=
\begin{bmatrix}e_{i_1,1}&e_{i_1,2}&e_{i_1,3}\\e_{i_2,1}&e_{i_2,2}&e_{i_2,3}\\e_{5,1}&e_{5,2}&e_{5,3}\\\end{bmatrix}=
\begin{bmatrix}
\hat e_{\hat i_1,1}&\hat e_{\hat i_1,2}& e_{i_1,3}\\
\hat e_{\hat i_2,1}&\hat e_{\hat i_2,2}&e_{i_2,3}\\e_{5,1}& e_{5,2}&e_{5,3}
\end{bmatrix}.
$$
We have
$$
E(I,J) \begin{bmatrix}I_2&-g\\0&1\end{bmatrix}=\begin{bmatrix}
\hat e_{\hat i_1,1}&\hat e_{\hat i_1,2}& 0\\
\hat e_{\hat i_2,1}&\hat e_{\hat i_2,2}&0\\e_{5,1}& e_{5,2}&e'_{5,3}
\end{bmatrix},
$$
where $e'_{5,3}=e_{5,3}-\begin{bmatrix}e_{5,1}&e_{5,2}\end{bmatrix}g\neq 0$;
hence $\det(E(I,J))=e'_{5,3}\det(\widehat E(\hat I,\hat J))\neq 0$.
\end{itemize}
}
\end{example}

{\bf Proof of Lemma \ref{lemmamingID}.}
If $k=r$, the result is trivial. Indeed, $Q_{r,r}=\{[1,\dots, r]\}$; therefore, $Z=Z^*=[1,\dots, r]$, and for  $I=J=Z$ we have $\det(E(I,J))=\det(E)\neq 0$.

If $k=1$, then $Z=[z_1]$ and $Z^*=[r-z_1+1]$. Since $\rank(E)=r$, we have $\rank(E(1:r-z_1+1,1:z_1)\geq 1$, therefore, there exist $i\leq r-z_1+1$ and $j\leq z_1$ such that $e_{i,j}\neq 0$. This means that Lemma \ref{lemmamingID} holds for $k=1$. Thus, we can assume from now on that $k\geq 2$ and, consequently, $r\geq 2$.

We  prove the lemma for $r\geq 2$  by induction on $r$.

If $r=2$,  then $k=r=2$ and, as mentioned above, the lemma is trivially satisfied.

Let $r-1 \geq 2$. Assume that the lemma is satisfied  for  $\widehat E\in{\mathcal R}^{(r-1)\times (r-1)}$ with $\rank (\widehat E) = r-1$,
and let $  E\in{\mathcal R}^{r\times r}$ with $\rank (E) = r$.

Let $X=  E(1:(r-1), :)$, i.e., $X\in{\mathcal R}^{(r-1)\times r}$ is the matrix formed by the $r-1$ first rows of
$  E$, and let $c_1, \dots, c_r\in {\mathcal R}^{(r-1)\times 1}$ be the columns of $X$. As $\rank (E)=r$, we have  $\rank(X)=r-1$.

If $c_1=0$, take $u=1$. Otherwise, let $c_u$ be the first column in $X$ that depends linearly on the preceding columns, i.e.,
$$
u=\min\{i: \rank\left(\begin{bmatrix}c_1&\cdots &c_i\end{bmatrix}\right)=i-1\}.
$$
Hence, $1\leq u \leq r$, and if $u>1$
there  exists $g\in {\mathcal F}^{(u-1)\times 1}$ such that
\begin{equation}\label{eqgc98}
\begin{bmatrix}c_1&\cdots &c_{u-1}\end{bmatrix}g=c_u.\end{equation}
As $\rank \left(\begin{bmatrix}c_1&\cdots &c_u\\e_{r,1}&\cdots &e_{r,u}\end{bmatrix}\right)=\rank (E(:, 1:u))=u$, we have
\begin{equation}\label{eqgne}
\begin{array}{ll}
0\neq e_{r,u},&\text{if } u=1,
\\
\begin{bmatrix}e_{r,1}&\dots &e_{r,u-1}\end{bmatrix}g\neq e_{r,u},&\text{if } u>1.
\end{array}
\end{equation}

Let $\widehat E=[\hat e_{i,j}]_{1\leq i,j\leq r-1} =\begin{bmatrix}c_1&\dots &c_{u-1}&c_{u+1}&\dots &c_{r} \end{bmatrix}\in{\mathcal R}^{(r-1)\times (r-1)}$, i.e., $\widehat E$ is the
submatrix obtained  from $E$  by deleting the last row and the column $u$.
 Then $\rank(\widehat E)=r-1$, and
$$\begin{array}{rll}
  \hat e_{i,j}=&e_{i,j}, & 1\leq i \leq r-1, \ 1\leq j \leq u-1,
  \\
  \hat e_{i,j}=&e_{i,j+1}, & 1\leq i \leq r-1, \ u\leq j \leq r-1.
  \end{array}
$$
Let $k \in \{2, \dots, r\}$ and $Z=[z_1, \dots, z_k]\in  Q_{k,r}$.
Set  $z_0=0$. Notice that $z_i\geq i$, $0\leq i \leq k$, and let $w=\max\{i\geq 0: z_i=i\}$. Observe that $0\leq w\leq k$.
\begin{itemize}

\item If $w< u$, take $\hat z_0=0$ and  $\hat Z=[\hat z_1, \dots, \hat z_k]$, where
$$
\begin{array}{l}
\hat z_i=z_i=i, \quad 1\leq i \leq w,\\
\hat z_i=z_i-1, \quad w+1\leq i \leq k.
\end{array}
$$
If $w=k$, then $ 1=\hat z_1< \hat z_k=k=w\leq u-1\leq r-1$.
If $w=0$, then $0<z_1-1=\hat z_1< \dots < \hat z_k=z_k-1\leq r-1$.
If $0<w<k$, then $0<1=z_1=\hat z_1\leq w=\hat z_{w}=w+1-1<z_{w+1}-1=\hat z_{w+1}\leq \dots \leq \hat z_k=z_k-1\leq r-1$. Hence, $\hat Z\in Q_{k, r-1}$.
Let $Z^*=[z^*_1, \dots, z^*_k]$ and  $\hat Z^*=[\hat z^*_1, \dots, \hat z^*_k]$. Then
$$
\begin{array}{ll}
\hat z^*_i=r-\hat z_{k-i+1}=r- z_{k-i+1}+1=z_{i}^*, & 1\leq i \leq k-w,\\
\hat z^*_i=r-\hat z_{k-i+1}=r- z_{k-i+1}=z_{i}^*-1, & k-w+1\leq i \leq k.
\end{array}
$$
 By the induction hypothesis,
there exist
 $\hat I=[\hat i_1, \dots, \hat i_k], \hat J=[\hat j_1, \dots, \hat j_k]\in Q_{k,r-1}\subset  Q_{k, r}$ such that $\hat I\leq  \hat Z^*$, $\hat J\leq \hat Z$,
and $\det (\widehat E(\hat I, \hat J))\neq 0$.
Taking
$I=\hat I$, we have $I\in Q_{k, r}$ and $I\leq \hat Z^*\leq Z^*$. Define  $J=[j_1, \dots, j_k]$ as
$$
\begin{array}{ll}
j_i=\hat j_i, & \mbox{ if } \hat j_i<u,\\
j_i=\hat j_i+1, &\mbox{ if } \hat j_i\geq u.
\end{array}
$$
Then,  $1\leq \hat j_1\leq j_1<\dots <j_k\leq \hat j_k+1\leq r$, hence  $J\in Q_{k, r}$. By the definition of $\widehat E$, $E(I,J)=\widehat E(\hat I,\hat J)$, therefore $\det (E(I,J))\neq 0$.
Moreover, if $\hat j_i<u$, then $j_i\leq \hat z_i\leq z_i$.
If $\hat j_i\geq u$, then $\hat z_i\geq \hat j_i \geq u>w=\hat z_w$; hence, $i\geq w+1$ and $j_i=\hat j_i+1\leq \hat z_i+1=z_i$.
 It means that $J\leq Z$.
\item If $w\geq u$, then $1\leq u \leq w \leq k$.
  Take $\hat z_0=0$ and define $\hat Z=[\hat z_1, \dots, \hat z_{k-1}]$ as
$$
\begin{array}{ll}
\hat z_i=z_i=i, & 1\leq i \leq u-1,\\
\hat z_i=z_{i+1}-1, & u\leq i \leq k-1.
\end{array}
$$
If $u=k$, then $1=\hat z_1\leq \hat z_{u-1}=\hat z_{k-1}=u-1\leq r-1$.
If $1<u<k$, then $1=\hat z_1\leq \hat z_{u-1}=u-1=z_{u-1}<z_u\leq z_{u+1}-1=\hat z_u\leq \hat z_{k-1}=z_{k}-1\leq r-1$.
If $u=1$, then $0<\hat z_1=z_2-1\leq \hat z_{k-1}=z_{k}-1\leq r-1$.
Therefore, $\hat Z\in Q_{k-1, r-1}$ .
Let $Z^*=[z^*_1, \dots, z^*_k]$ and  $\hat Z^*=[\hat z^*_1, \dots, \hat z^*_{k-1}]$. Then,
$$
\begin{array}{ll}
\hat z^*_i=r-\hat z_{k-i}=r- z_{k-i+1}+1=z^*_i, & 1\leq i \leq k-u,
\\
\hat z^*_i=r-\hat z_{k-i}=
r-k+i=r+1-z_{k-i+1}=z^*_i,
& k-u+1\leq i \leq k-1.
\end{array}
$$
By the induction hypothesis, there exist
$\hat I=[\hat i_1, \dots, \hat i_{k-1}]$, $\hat J=[\hat j_1, \dots, \hat j_{k-1}]\in Q_{k-1,r-1}$ such that $\hat I\leq  \hat Z^*$, $\hat J\leq \hat Z$,
and $\det (\widehat E(\hat I, \hat J))\neq 0$. Note that $\hat j_i = \hat z_i = i$, for $1\leq i \leq u-1$.
Define $I=[i_1, \dots, i_k]$ as
$$
\begin{array}{l}
i_h=\hat i_h, \quad 1\leq h \leq k-1,\\
i_k=r.
\end{array}
$$
Then, $I\in Q_{k, r}$ and $i_h=\hat i_h\leq \hat z^*_h =z^*_h$ for $1\leq h \leq k-1$. Since $1\leq w$, we have $z_1=1$ and $z^*_k=r$; hence $i_k = z^*_k$ and $I\leq Z^*$.

Define also  $J=[j_1, \dots, j_k]$ as
$$
\begin{array}{ll}
j_i=i, &1\leq i \leq u,\\
j_i=\hat j_{i-1}+1, &u+1\leq i \leq k.
\end{array}
$$
 As $j_u=u\leq \hat j_u <\hat j_u+1 = j_{u+1}$, we have $1=j_1<\dots <j_k\leq r$.
 Therefore, $J\in Q_{k, r}$
and
$$ E(I,J)=\begin{bmatrix}
    e_{i_1, 1}& \dots &e_{i_1, u}&e_{i_1, j_{u+1}}&\dots &e_{i_1, j_k}\\
    \vdots &\ddots &\vdots&\vdots&\ddots &\vdots\\
    e_{i_{k-1}, 1}& \dots &e_{i_{k-1}, u}&e_{i_{k-1}, j_{u+1}}&\dots &e_{i_{k-1}, j_k}\\ e_{r, 1}& \dots &e_{r, u}&e_{r, j_{u+1}}&\dots &e_{r, j_k}\\
  \end{bmatrix}.$$
Thus, taking into account \eqref{eqgc98},
  \begin{align*}
& E(I,J)\begin{bmatrix}
    I_{u-1}&-g&0\\
    0&1&0\\
    0&0&I_{k-u}
    \end{bmatrix} \\ & \qquad =\begin{bmatrix}
    e_{i_1, 1}& \dots &e_{i_1, u-1}&0&e_{i_1, j_{u+1}}&\dots &e_{i_1, j_k}\\
    \vdots &\ddots &\vdots&\vdots&\vdots&\ddots &\vdots\\
    e_{i_{k-1}, 1}& \dots &e_{i_{k-1}, u-1}&0&e_{i_{k-1}, j_{u+1}}&\dots &e_{i_{k-1}, j_k}\\ e_{r, 1}& \dots &e_{r, u-1}&e'_{r, u}&e_{r, j_{u+1}}&\dots &e_{r, j_k}\\
  \end{bmatrix},
  \end{align*}
  where $e'_{r, u}=e_{r,u}-\begin{bmatrix}e_{r,1}&\dots &e_{r,u-1}\end{bmatrix}g$
   if $u>1$, and $e'_{r,u}=e_{r,u}$ if $u=1$.
From (\ref{eqgne}) we obtain  $e'_{r,u}\neq 0$.
Moreover,
$$
\begin{array}{rl}
&\begin{bmatrix}
    e_{i_1, 1}& \dots &e_{i_1, u-1}&e_{i_1, j_{u+1}}&\dots &e_{i_1, j_k}\\
    \vdots &\ddots &\vdots&\vdots&\ddots &\vdots\\
    e_{i_{k-1}, 1}& \dots &e_{i_{k-1}, u-1}&e_{i_{k-1}, j_{u+1}}&\dots &e_{i_{k-1}, j_k}
  \end{bmatrix}
\vspace{7pt}\\ \vspace{7pt}
=&\begin{bmatrix}
    e_{\hat i_1, \hat j_1}& \dots &e_{\hat i_1, \hat j_{u-1}}&e_{\hat i_1, \hat j_{u}+1}&\dots &e_{\hat i_1, \hat j_{k-1}+1}\\
    \vdots &\ddots &\vdots&\vdots&\ddots &\vdots\\
    e_{\hat i_{k-1}, \hat j_1}& \dots &e_{\hat i_{k-1},  \hat j_{u-1}}&e_{\hat i_{k-1}, \hat j_{u}+1}&\dots &e_{\hat i_{k-1}, \hat j_{k-1}+1}
  \end{bmatrix}
\vspace{7pt}\\ \vspace{7pt}
=&\begin{bmatrix}
    \hat e_{\hat i_1, \hat j_1}& \dots &\hat e_{\hat i_1, \hat j_{u-1}}&\hat e_{\hat i_1, \hat j_{u}}&\dots &\hat e_{\hat i_1, \hat j_{k-1}}\\
    \vdots &\ddots &\vdots&\vdots&\ddots &\vdots\\
    \hat e_{\hat i_{k-1}, \hat j_1}& \dots &\hat e_{\hat i_{k-1},  \hat j_{u-1}}&\hat e_{\hat i_{k-1}, \hat j_{u}}&\dots &\hat e_{\hat i_{k-1}, \hat j_{k-1}}
  \end{bmatrix}
\\
=&
\widehat E(\hat I,\hat J).
\end{array}
$$
Hence
$\det(E(I,J))=(-1)^{k+u}e'_{r,u}
\det(\widehat E(\hat I,\hat J))\neq 0$.
Moreover, $j_i=i= z_i$ for $1\leq i \leq u$, and $j_i=\hat j_{i-1}+1\leq \hat z_{i-1}+1=z_i$ for $u+1\leq i \leq k$.
Thus, $J\leq Z$.
  \end{itemize}
\hfill $\Box$

\end{appendices}

\bibliographystyle{acm}

\begin{thebibliography}{10}

\bibitem{AmBaMaRo23}
{\sc Amparan, A., {Baraga\~na}, I., Marcaida, S., and Roca, A.}
\newblock Row or column completion of polynomial matrices of given degree.
\newblock {\em SIAM J. Matrix Anal. Appl. 45}, 1 (2024), 478--503.

\bibitem{AmBaMaRo25}
{\sc Amparan, A., {Baraga\~na}, I., Marcaida, S., and Roca, A.}
\newblock Row completion of polynomial and rational matrices.
\newblock {\em Submitted. Available in https://arxiv.org/pdf/2504.10303\/}
  (2025).

\bibitem{AmMaZa15}
{\sc Amparan, A., Marcaida, S., and Zaballa, I.}
\newblock Finite and infinite structures of rational matrices: a local
  approach.
\newblock {\em Electron. J. Linear Algebra 30\/} (2015), 196--226.

\bibitem{amparan2024parametrizing}
{\sc Amparan, A., Marcaida, S., and Zaballa, I.}
\newblock Parametrizing spectral filters for nonsingular polynomial matrices.
\newblock {\em Linear and Multilinear Algebra\/} (published online in 2024),
  1--45.

\bibitem{anguasetal2019}
{\sc Anguas, L.~M., Dopico, F.~M., Hollister, R., and Mackey, D.~S.}
\newblock Van {D}ooren's index sum theorem and rational matrices with
  prescribed structural data.
\newblock {\em SIAM J. Matrix Anal. Appl. 40}, 2 (2019), 720--738.

\bibitem{BaZa02}
{\sc {Baraga\~na}, I., and Zaballa, I.}
\newblock Feedback invariants of restrictions and quotients: series connected
  systems.
\newblock {\em Linear Algebra Appl. 351--352\/} (2002), 69--89.

\bibitem{batzke-mehl-2014}
{\sc Batzke, L., and Mehl, C.}
\newblock On the inverse eigenvalue problem for {$T$}-alternating and
  {$T$}-palindromic matrix polynomials.
\newblock {\em Linear Algebra Appl. 452\/} (2014), 172--191.

\bibitem{nlevp-collection}
{\sc Betcke, T., Higham, N.~J., Mehrmann, V., Schr\"oder, C., and Tisseur, F.}
\newblock N{LEVP}: a collection of nonlinear eigenvalue problems.
\newblock {\em ACM Trans. Math. Software 39}, 2 (2013), Art. 7, 1--28.

\bibitem{das-bora}
{\sc Das, B., and Bora, S.}
\newblock Nearest rank deficient matrix polynomials.
\newblock {\em Linear Algebra Appl. 674\/} (2023), 304--350.

\bibitem{DeTDmyDop2020}
{\sc De~Ter\'an, F., Dmytryshyn, A., and Dopico, F.~M.}
\newblock Generic symmetric matrix polynomials with bounded rank and fixed odd
  grade.
\newblock {\em SIAM J. Matrix Anal. Appl. 41}, 3 (2020), 1033--1058.

\bibitem{DeDoMa14}
{\sc {De Ter\'an}, F., Dopico, F.~M., and Mackey, D.~S.}
\newblock Spectral equivalence of matrix polynomials and the index sum theorem.
\newblock {\em Linear Algebra Appl. 459\/} (2014), 264--333.

\bibitem{quadratic-palin}
{\sc De~Ter\'an, F., Dopico, F.~M., Mackey, D.~S., and Perovi\'c, V.}
\newblock Quadratic realizability of palindromic matrix polynomials.
\newblock {\em Linear Algebra Appl. 567\/} (2019), 202--262.

\bibitem{DeDoMaVa16}
{\sc {De Ter\'an}, F., Dopico, F.~M., Mackey, D.~S., and {Van Dooren}, P.}
\newblock Polynomial zigzag matrices, dual minimal bases, and the realization
  of completely singular polynomials.
\newblock {\em Linear Algebra Appl. 488\/} (2016), 460--504.

\bibitem{DeDoVa15}
{\sc {De Ter\'an}, F., Dopico, F.~M., and {Van Dooren}, P.}
\newblock Matrix polynomials with completely prescribed eigenstructure.
\newblock {\em SIAM J. Matrix Anal. Appl. 36}, 1 (2015), 302--328.

\bibitem{DmDoVa23}
{\sc Dmytryshyn, A., Dopico, F.~M., and Van~Dooren, P.}
\newblock Minimal rank factorizations of polynomial matrices.
\newblock {\em Linear Algebra Appl. in press
  (https://doi.org/10.1016/j.laa.2025.01.009)\/} (2025).

\bibitem{dmytryshyn-geom}
{\sc Dmytryshyn, A., Johansson, S., K{\aa}gstr\"om, B., and Van~Dooren, P.}
\newblock Geometry of matrix polynomial spaces.
\newblock {\em Found. Comput. Math. 20\/} (2020), 423--450.

\bibitem{local-lin}
{\sc Dopico, F.~M., Marcaida, S., Quintana, M.~C., and Van~Dooren, P.}
\newblock Local linearizations of rational matrices with application to
  rational approximations of nonlinear eigenvalue problems.
\newblock {\em Linear Algebra Appl. 604\/} (2020), 441--475.

\bibitem{Fo75}
{\sc Forney, G.~D.}
\newblock Minimal bases of rational vector spaces, with applications to
  multivariable linear systems.
\newblock {\em SIAM J. Control 13}, 3 (1975), 493--520.

\bibitem{FuWi79}
{\sc Fuhrmann, P., and Willems, J.}
\newblock Factorization indices at infinity for rational matrix functions.
\newblock {\em Integral Equations Operator Theory 2}, 3 (1979), 287--301.

\bibitem{GLR-SIAM-2009}
{\sc Gohberg, I., Lancaster, P., and Rodman, L.}
\newblock {\em Matrix {P}olynomials}, vol.~58 of {\em Classics in Applied
  Mathematics}.
\newblock Society for Industrial and Applied Mathematics (SIAM), Philadelphia,
  PA, 2009.
\newblock Reprint of the 1982 original.

\bibitem{guttel-tisseur-nlep}
{\sc G\"uttel, S., and Tisseur, F.}
\newblock The nonlinear eigenvalue problem.
\newblock {\em Acta Numer. 26\/} (2017), 1--94.

\bibitem{Kail80}
{\sc Kailath, T.}
\newblock {\em Linear Systems}.
\newblock Prentice Hall, New Jersey, 1980.

\bibitem{LaTi85}
{\sc Lancaster, P., and Tismenetsky, M.}
\newblock {\em The Theory of Matrices}, second~ed.
\newblock Academic Press, Inc., Orlando, FL, 1985.

\bibitem{MMMM-good}
{\sc Mackey, D.~S., Mackey, N., Mehl, C., and Mehrmann, V.}
\newblock Structured polynomial eigenvalue problems: good vibrations from good
  linearizations.
\newblock {\em SIAM J. Matrix Anal. Appl. 28}, 4 (2006), 1029--1051.

\bibitem{MMMM-lin-class}
{\sc Mackey, D.~S., Mackey, N., Mehl, C., and Mehrmann, V.}
\newblock Vector spaces of linearizations for matrix polynomials.
\newblock {\em SIAM J. Matrix Anal. Appl. 28}, 4 (2006), 971--1004.

\bibitem{MMMM15}
{\sc Mackey, D.~S., Mackey, N., Mehl, C., and Mehrmann, V.}
\newblock Möbius transformations of matrix polynomials.
\newblock {\em Linear Algebra Appl. 470\/} (2015), 120--184.

\bibitem{Sa80}
{\sc {Marques de S\`a}, E.}
\newblock On the diagonals of integral matrices.
\newblock {\em Czechoslovak Mathematical Journal 30}, 2 (1980), 207--212.

\bibitem{MOA11}
{\sc Marshall, A.~W., Olkin, I., and Arnold, B.~C.}
\newblock {\em Inequalities: Theory of Majorization and its Applications},
  second~ed.
\newblock Springer, New York, 2011.

\bibitem{root-vectors-vanni-VD}
{\sc Noferini, V., and Van~Dooren, P.}
\newblock Root vectors of polynomial and rational matrices: theory and
  computation.
\newblock {\em Linear Algebra Appl. 656\/} (2023), 510--540.

\bibitem{perovic-mackey}
{\sc Perovi\'c, V., and Mackey, D.~S.}
\newblock Quadratic realizability of palindromic matrix polynomials: the real
  case.
\newblock {\em Linear Multilinear Algebra 71}, 5 (2023), 797--841.

\bibitem{Rose70}
{\sc Rosenbrock, H.~H.}
\newblock {\em State-space and Multivariable Theory}.
\newblock Thomas Nelson and Sons, London, 1970.

\bibitem{strangfund}
{\sc Strang, G.}
\newblock The fundamental theorem of linear algebra.
\newblock {\em Amer. Math. Monthly 100}, 9 (1993), 848--855.

\bibitem{taslaman-tiss-zab}
{\sc Taslaman, L., Tisseur, F., and Zaballa, I.}
\newblock Triangularizing matrix polynomials.
\newblock {\em Linear Algebra Appl. 439}, 7 (2013), 1679--1699.

\bibitem{tisseur-zaballa}
{\sc Tisseur, F., and Zaballa, I.}
\newblock Triangularizing quadratic matrix polynomials.
\newblock {\em SIAM J. Matrix Anal. Appl. 34}, 2 (2013), 312--337.

\bibitem{vandooren-thesis}
{\sc Van~Dooren, P.}
\newblock {\em The Generalized Eigenstructure Problem: Applications in Linear
  System Theory}.
\newblock Ph.D. dissertation, Kath. Univ. Leuven, Belgium, 1979.

\bibitem{vandooren-iee81}
{\sc Van~Dooren, P.}
\newblock The generalized eigenstructure problem in linear system theory.
\newblock {\em IEEE Trans. Automat. Control 26}, 1 (1981), 111--129.

\bibitem{vd-dewilde83}
{\sc Van~Dooren, P., and Dewilde, P.}
\newblock The eigenstructure of an arbitrary polynomial matrix: computational
  aspects.
\newblock {\em Linear Algebra Appl. 50\/} (1983), 545--579.

\bibitem{Vard91}
{\sc Vardulakis, A. I.~G.}
\newblock {\em Linear Multivariable Control}.
\newblock John Wiley and Sons, New York, 1991.

\bibitem{vergheseetal}
{\sc Verghese, G., Van~Dooren, P., and Kailath, T.}
\newblock Properties of the system matrix of a generalized state-space system.
\newblock {\em Internat. J. Control 30}, 2 (1979), 235--243.

\bibitem{Wo74}
{\sc Wolovich, W.}
\newblock {\em Linear Multivariable Systems}.
\newblock Springer, New York-Heidelberg, 1974.

\bibitem{Za97}
{\sc Zaballa, I.}
\newblock Controllability and {H}ermite indices of matrix pairs.
\newblock {\em Int. J. Control 68}, 1 (1997), 61--86.

\end{thebibliography}

\end{document}